\documentclass[preprint,12pt]{elsarticle}

\usepackage{epstopdf}
\ifpdf
  \DeclareGraphicsExtensions{.eps,.pdf,.png,.jpg}
\else
  \DeclareGraphicsExtensions{.eps}
\fi

\graphicspath{{figs/}}

\usepackage{amssymb}
\usepackage{amsthm}

\newtheorem{theorem}{Theorem}
\newtheorem{corollary}{Corollary}

\usepackage{amsmath}
\usepackage{verbatim}

\journal{Applied Numerical Mathematics}

\begin{document}
%% \linenumbers
 
\begin{frontmatter}

%% Title, authors and addresses

%% use the tnoteref command within \title for footnotes;
%% use the tnotetext command for theassociated footnote;
%% use the fnref command within \author or \address for footnotes;
%% use the fntext command for theassociated footnote;
%% use the corref command within \author for corresponding author footnotes;
%% use the cortext command for theassociated footnote;
%% use the ead command for the email address,
%% and the form \ead[url] for the home page:
%% \title{Title\tnoteref{label1}}
%% \tnotetext[label1]{}
%% \author{Name\corref{cor1}\fnref{label2}}
%% \ead{email address}
%% \ead[url]{home page}
%% \fntext[label2]{}
%% \cortext[cor1]{}
%% \address{Address\fnref{label3}}
%% \fntext[label3]{}

\title{Spectral computation of low probability tails for the homogeneous Boltzmann equation}

%% use optional labels to link authors explicitly to addresses:
%% \author[label1,label2]{}
%% \address[label1]{}
%% \address[label2]{}

\author[UTDMath]{John Zweck\corref{cor1}}
\ead{zweck@utdallas.edu}

\author[UTDMath]{Yanping Chen}
\ead{yanpingchen123@yahoo.com}

\author[UTDPhysics]{Matthew J. Goeckner}
\ead{goeckner@utdallas.edu}

\author[KSMath]{Yannan Shen}
\ead{yshen@ku.edu}

\cortext[cor1]{Corresponding author}

\address[UTDMath]{Department of Mathematical Sciences, The University of Texas at Dallas, Richardson, TX 75080, USA}

\address[UTDPhysics]{Department of Physics, The University of Texas at Dallas, Richardson, TX 75080, USA}

\address[KSMath]{Department of Mathematics, University of Kansas, Lawrence, KS 66045, USA}

\begin{abstract}

%% Text of abstract
% REQUIRED
We apply the spectral-Lagrangian method of Gamba and Tharkabhushanam
for solving the homogeneous Boltzmann equation 
to compute the low probability tails of the
velocity distribution function, $f$,  of a particle species. 
This method is based on a truncation, $Q^{\operatorname{tr}}(f,f)$,
of the Boltzmann collision operator, $Q(f,f)$,  whose Fourier transform 
is given by  a weighted convolution. The truncated collision operator
models the situation in which two colliding particles ignore each other
if their relative speed exceeds a threshold, $g_{\text{tr}}$.
We demonstrate that the choice of truncation 
parameter plays a
critical role in the accuracy of the numerical computation of $Q$.
Significantly, if $g_{\text{tr}}$ is too large, then accurate numerical computation of the 
weighted convolution integral  is not feasible,
since the decay rate and degree of oscillation of the convolution weighting function 
both increase as $g_{\text{tr}}$ increases.
We derive an upper bound on the pointwise
error between $Q$ and $Q^{\text{tr}}$,
assuming that both operators are computed exactly. 
This bound provides some additional theoretical justification for the spectral-Lagrangian
method, and can be used to guide the choice of $g_{\text{tr}}$
in numerical computations. We  then
demonstrate  how to choose $g_{\text{tr}}$ and the numerical 
discretization parameters 
so that the computation of the truncated collision operator 
is a good approximation to $Q$ in the low probability tails.
Finally, for several different initial conditions, we 
demonstrate the feasibility of accurately computing the  time evolution
of the velocity pdf  down to probability density levels 
ranging from $10^{-5}$ to $10^{-9}$.
\end{abstract}

\begin{keyword}
%% keywords here, in the form: keyword \sep keyword
Boltzmann collision operator \sep spectral numerical method \sep  low-probability tails 
%% PACS codes here, in the form: \PACS code \sep code

%% MSC codes here, in the form: \MSC code \sep code
%% or \MSC[2008] code \sep code (2000 is the default)
\MSC[2020] 35Q20 \sep 35R09 \sep 82C40 \sep 82D10 \sep 65Z05
\end{keyword}

\end{frontmatter}

\section{Introduction}

The motivation for this work is to develop improved computational tools
for the simulation of low-probability, high-energy processes in non-equilibrium, 
low-temperature plasmas. Our interest is in kinetic models 
for the evolution of the velocity  probability density function (pdf) 
of each particle species in a plasma. Such models are based on the 
Boltzmann equation which governs both the transport of, and collisions
between, particles. In plasma systems, gas-phase chemistry and surface kinetics
 are largely driven by collision processes between high-energy electrons in the plasma and molecules in the gas phase~\cite{Gustavsson:2004}.  Reaction rates in the gas phase are determined by the overlap between the electron velocity pdf
 and the electron-impact cross sections of the various species. 
Accurate calculation of the low-probability tails of the electron velocity pdf 
is therefore critical. 
If the plasma is in thermal equilibrium, the 
electron velocity pdf can often be assumed to be Maxwellian.
However, experimental results demonstrate that  the  Maxwellian assumption
is often invalid~\cite{allen1978applicability,Dicarlo:1989,sheridan1998electron,Sozzi:2008,tan1973langmuir}, especially for pulsed plasmas where 
the velocity pdf may depend strongly on both spatial position and on time~\cite{poulose2017driving}.

The Direct Simulation Monte Carlo method (DSMC) is often used to numerically model collision processes in inhomogeneous (position-dependent) systems, and in systems that are not in thermal equilibrium. This method was initially developed by 
Bird~\cite{Bird:1994} and Nanbu~\cite{Nanbu:1983}. Wagner proved
that solutions obtained using the DSMC method converge to the solution
of the Boltzmann equation~\cite{Wagner:1992}. 
More recently, Rjasanow, Gamba, and Wagner modified the DSMC method  to compute the low-probability tails of steady state solutions~\cite{Gamba:2005,Rjasanow:2005}. 
Although they have proved effective in many situations, the statistical 
uncertainties in these methods can be challenging to resolve for systems that are not in
thermal equilibrium~\cite{GambaThark228p2012}.

Rather than attempting to model a realistic plasma system, in this paper we
focus on the narrower goal of computing the 
velocity pdf, $f$, of a particle species (such as the electrons)
down into the low-probability tails
under the assumption that $f=f(t,\boldsymbol v)$ satisfies an initial-value
problem for the homogeneous Boltzmann equation,
\begin{equation}
\frac{\partial f}{\partial t}\,\,=\,\, Q(f,f).
\label{eq:HBEIntro}
\end{equation}
Here, the Boltzmann collision operator, $Q$, is a bilinear integral operator
that  is defined in terms of a kernel that models a binary collision process.
Although it omits much of the physics, this computation is nevertheless
challenging because for each time, $t$, and each 
point in a 3-dimensional space of velocities, $\boldsymbol v$, 
the evaluation of $Q(f,f)(t,\boldsymbol v)$ involves 
the computation of a 5-dimensional integral over a space of velocities and angular directions, resulting in a computational cost of order $\mathcal O(N^8)$.

Over the last two decades 
there have been several major advances that have enabled 
more efficient  computation of the Boltzmann collision operator.
An important class of deterministic methods are the spectral methods
which include the Fourier-Galerkin methods of Pareschi and his
collaborators~\cite{Mouhot:2006,pareschi1996fourier,Pareschi:2000}, 
the spectral-Lagrangian methods of Gamba and her group~\cite{gamba2014conservative,gamba2017fast,GambaThark228p2012,GambaThark28p430,haack2012high,munafo2014spectral}, and the more recent Petrov-Galerkin method
of Gamba and Rjasanow~\cite{gamba2018galerkin}. 

With the Fourier-Galerkin method of Pareschi and Russo~\cite{Pareschi:2000}, the
velocity pdf is assumed to be compactly 
supported and is approximated by a finite Fourier series. The Boltzmann collision operator then takes the form of a  weighted discrete
convolution operator where the weights are given in terms of the collision kernel. 
The resulting numerical scheme has a computational cost of $\mathcal O(N^6)$,
where $N$ is the number of discretization points in each velocity
dimension, which represents a substantial improvement over the 
$\mathcal O(N^8)$ cost of direct numerical integration
of the collision operator. 
Moreover, the method is spectrally accurate  and conserves mass. 
However, due to the use of a Fourier series representation, positivity of the
solution is not guaranteed, and non-physical high energy collisions are
incorporated into the model due to the periodization and truncation of the
velocity pdf and the collision operator.
%%
%In subsequent work, Mouhot and Pareschi~\cite{Mouhot:2006} 
%developed fast algorithms for $Q$ that rely on a Carelman-like
%representation to approximate the collision operator as a linear combination of pure convolutions.
%These methods can be applied to particular classes of collision processes
%and have been shown to be competitive with Monte Carlo simulations 
%even in the spatially inhomogeneous case~\cite{Filbet:2006}.
%For example, for hard-sphere collisions in three velocity 
%dimensions the resulting algorithm 
%is $\mathcal O(A N^3 \log N)$, where $A \ll N^3$. 
%%
Building on this  approach, 
Gamba et al.~\cite{gamba2017fast} developed a 
$\mathcal O(M N^4 \log N)$  algorithm with $M \ll N^2$, valid for arbitrary collision kernels, in which a pure convolution structure is achieved by numerical quadrature
of the integral defining the convolution weighting function. 
Other advances along these
lines include a method of Fonn et al.~\cite{fonn2014hyperbolic} that operates on a sparse  set of Fourier modes, and a method of Cai et al.~\cite{cai2018entropic} that  preserves positivity at Fourier collocation points and satisfies the H-theorem. 

Compared to the Fourier-Galerkin methods,
the spectral-Lagrangian method of Gamba and Tharkabhushanam~\cite{GambaThark228p2012}
has the advantage that it provides a general framework for
arbitrary collision kernels with either elastic or
inelastic binary interactions, does not require 
periodization of $f$, and enforces conservation of moments
through solution of an auxiliary constrained optimization problem.  
The method is based on a formula for the Fourier transform of the
collision operator in the form of a weighted convolution, 
\begin{equation} \label{QhatIntro}
\widehat{Q}(f,f)(\boldsymbol{\zeta}) = (2\pi)^{-3/2}
\int_{\mathbb{R}^3} \widehat{f}(\boldsymbol{\zeta} - \boldsymbol{\xi}) 
\widehat{f}(\boldsymbol{\xi}) \widehat{G}(\boldsymbol{\xi},\boldsymbol{\zeta}) \, d\boldsymbol{\xi},
\end{equation}
where  $\widehat{G}$ is a  convolution weighting function
that can be precomputed.
The computational cost of the method is therefore the cost of
numerically computing the integrals \eqref{QhatIntro} for all $\boldsymbol\zeta \in \mathbb R^3$, which is $\mathcal O(N^6)$.
Analogous to~\cite{gamba2017fast}, 
Gamba et al.~\cite{gamba2014fast}
obtained an approximate formula for $\widehat{G}$
which enables~\eqref{QhatIntro} to be expressed as a pure convolution that can 
be sped up using the fast Fourier transform to yield a 
$\mathcal O(M N^4 \log N)$ algorithm with $M \ll N^2$.

Alonso, Gamba, and Tharkabhushanam~\cite{alonso2018convergence}
analyzed the accuracy and consistency of the spectral-Lagrangian method.
They restricted  $f$ and $Q(f,f)$  to a finite rectangular domain, 
$\Omega_L\subset \mathbb R^3$, of side-length, $L$, in velocity space
and then orthogonally projected onto 
an $N$-dimensional Fourier series basis yielding the initial value problem,
\begin{equation}
\frac{\partial h}{\partial t} \,\,=\,\, \Pi^N Q(h,h), \qquad \text{in } (0,T]\times \Omega_L,
\label{HBEProj}
\end{equation}
with $h(0,\boldsymbol v) = \Pi^N f(0,\boldsymbol v) $. (Here $\Pi^N$ is the projection operator.)
To enforce conservation of 
mass, momentum, and energy (for elastic collisions), 
they used the method of Lagrange multipliers to replace the the right-hand side
of~\eqref{HBEProj}  by the $L^2(\Omega_L)$-closest function  to $\Pi^N Q(h,h)$
with zero mass, momentum, and energy. They then proved 
that for a large class of initial data, one can choose the size of the truncated
domain, $\Omega_L$, the number of Fourier modes, $N$, and the final simulation time, $T$,  so that the solution, $h$, agrees with the
equilibrium Maxwellian distribution to within a desired tolerance in a suitable Sobolev norm.

The convolution weighting function,
$\widehat G(\boldsymbol{\xi},\boldsymbol{\zeta})$, in \eqref{QhatIntro}
is given  as the Fourier transform with respect to $\boldsymbol g$ of a 
kernel, $G(\boldsymbol{\xi},\boldsymbol g)$. To avoid the introduction of a 
divergent improper integral, 
the integral defining this Fourier transform 
must be taken over a finite ball, $|\boldsymbol g | \leq g_{\operatorname{tr}}$,
rather than over all of $\mathbb R^3$. 
Therefore, the spectral-Lagrangian method is based on an 
approximation, $Q^{\text{tr}}$, of $Q$, which we refer to as the truncated
collision operator. 
In their proof of  an existence and uniqueness
theorem for solutions of \eqref{eq:HBEIntro},
Cercignani et al.~\cite{cercignani2013mathematical} 
show that   $Q^{\text{tr}}$ converges weakly to $Q$ as  $g_{\text{tr}}\to\infty$.
Physically,  $Q^{\text{tr}}$ models the situation in which
two colliding particles ignore each other if their relative speed exceeds
the threshold, $g_{\text{tr}}$~\cite{cercignani2013mathematical}.
Pareschi and Russo~\cite{Pareschi:2000}, showed 
that if the velocity pdf has compact support 
in a ball of radius $R$ then $Q=Q^{\text{tr}}$ provided that $g_{\text{tr}} \geq 2R$.
They used this observation to avoid aliasing in their method to compute $Q$ using a 
Fourier series approximation of $f$.\footnote{We note that the Fourier series method
of Pareschi and Russo has  a different character to the spectral-Lagrangian method
of Gamba and Tharkabhushanam, which does not have to avoid the possibility of aliasing
effects.}
In a similar vein, Gamba and Tharkabhushanam showed that if the velocity pdf has compact support 
in the box, $\Omega_L$, then $Q=Q^{\text{tr}}$ provided that $g_{\text{tr}} \geq 2\sqrt{3}L$.
However, in their analysis of the method,  
Alonso et al.~\cite{alonso2018convergence}
assume that the function, $Q^{\text{tr}}(f,f)$, is computed
exactly from $f$, that is, they do not analyze the error in the numerical
computation of the integral~\eqref{QhatIntro} for the Fourier transform of the truncated collision operator.

The first goal of this paper is to demonstrate that
with the method of Gamba and Tharkabhushanam,
 the choice of the truncation 
parameter, $g_{\text{tr}}$,  plays a
critical role in the accuracy of the numerical computation of $Q$.
Clearly, if $g_{\text{tr}}$  is too small 
then $Q^{\text{tr}}$ will not be  a good approximation to $Q$. 
However, if $g_{\text{tr}}$ is too large then accurate numerical computation of the 
convolution integral~\eqref{QhatIntro} is not possible 
since the convolution weighting function, $\widehat{G}$,  is 
a slowly decaying oscillatory function of $\boldsymbol \xi$
whose degree of oscillation increases as $g_{\text{tr}}$ increases.
Indeed, with 
Gamba and Tharkabhushanam's theoretical choice of
$g_{\text{tr}} = 2\sqrt{3}L$, 
we show that the numerically computed
collision operator is a poor approximation. 
In unpublished work, Haack~\cite{Haack:2013} instead uses 
$g_{\text{tr}}= L$. 
However he provides no explanation for the smaller choice of $g_{\text{tr}}$.

Our second goal is to 
derive an upper bound on the pointwise
error between $Q$ and $Q^{\text{tr}}$,
assuming that both operators are computed exactly. 
This error estimate can be viewed as a generalization to velocity pdfs
without compact support of the
formula for  $g_{\text{tr}}$ obtained by  Gamba and Tharkabhushanam.  
In particular our estimate 
yields the following \emph{simple strategy} for  choosing the parameter, $g_{\operatorname{tr}}$,
in numerical computations of the low-probability tails. 
Specifically, to guarantee  that $Q^{\operatorname{tr}}$
is an accurate approximation to $Q$ at $\mathbf v$, we should choose
$g_{\operatorname{tr}}$ to be slightly larger than $|\mathbf v|$, \emph{irrespective of the
probability level at $\mathbf v$.}  
In particular, our error bound provides a theoretical justification for 
Haack's choice of $g_{\text{tr}}= L$. We obtained this  error bound 
without regard to any  truncation or discretization of the domain
of the velocity pdf or the collision operator. In particular, our  bound
has nothing to do with the need to avoid aliasing in the Fourier series
approximations of $f$ in methods such as that of Pareschi and Russo.
Indeed, we emphasize  that the method of Gamba and Tharkabhushanam
that we are analyzing here is based on a formula for 
$Q$ in terms of a continuous Fourier transform and so does not
require any periodization of the domain of $f$~\cite{GambaThark228p2012}.

Our  third goal is to perform a series of simulation studies that 
demonstrate  how to choose $g_{\text{tr}}$ and the numerical 
discretization parameters 
so that the numerical computation of the truncated collision operator 
is a good approximation to $Q$ in the low probability tails.
In particular, we demonstrate that when we use the 
simple strategy described above to select $g_{\text{tr}}$  to
guarantee  that $Q^{\operatorname{tr}}$  is close to $Q$ at $\mathbf v$, then 
it is often feasible to numerically compute the
generalized convolution integral  for $Q^{\operatorname{tr}}$ 
with sufficient accuracy at $\mathbf v$. 
Finally, 
for several different initial conditions, we use the selected values of 
 $g_{\text{tr}}$ to show that the time evolution
of the velocity pdf
can be computed accurately down to probability density levels 
ranging from $10^{-5}$ to $10^{-9}$.

In Section~\ref{Sec:SpecMethod}, we review the spectral-Lagrangian method of 
Gamba and Tharkabhushanam, and in Section~\ref{Sec:Prelim} we present
the results of preliminary numerical simulations  that demonstrate that the 
choice $g_{\text{tr}}$  plays a
critical role in the accuracy of the numerical computation of $Q$.
In Section~\ref{Sec:Tr}, we  derive the
bound on the relative error between $Q^{\operatorname{tr}}$ and $Q$.
In Section~\ref{Sec:Method}, we discuss some implementation details,
and in Section~\ref{Sec:Results} we present the results of our numerical simulations.
Finally, in Section~\ref{Sec:Conc} we make some conclusions.

\section{The spectral-Lagrangian method for the Boltzmann equation}\label{Sec:SpecMethod}

In this section, we review the spectral-Lagrangian method
for the homogeneous Boltzmann equation developed by
Gamba and Tharkabhushanam~\cite{GambaThark228p2012}.
This method reduces the computational cost of the collision operator
from $\mathcal O(N^8)$ to $\mathcal O(N^6)$, where $N$ is the number of discretization points in each velocity dimension. 

The homogeneous Boltzmann equation for the velocity probability density 
function (pdf), $f=f(t,\boldsymbol v)$, of particles of species due to elastic collisions with particles of the same species is given by
\begin{equation}
\frac{\partial f}{\partial t}\,\,=\,\, Q(f,f),
\label{eq:HBE}
\end{equation}
where the Boltzmann collision operator, $Q := 
Q(f,f)$,
 is given by
\begin{equation}  \label{elasticQ} %%\label{eq:Qtr}
Q(\boldsymbol{v})=\int\limits_{\mathbb{R}^3} 
 \int\limits_{S^2}
\left[ f(\boldsymbol{v}')f(\boldsymbol{w}')-f(\boldsymbol{v})f(\boldsymbol{v} +\boldsymbol{g}) \right] B\left(g,\frac{\boldsymbol{g} \cdot \Theta}{g}\right) d\Theta d\boldsymbol{g}.
\end{equation}
Here $\boldsymbol g = \boldsymbol w - \boldsymbol v$ is 
the relative pre-collisional velocity and $g=|\boldsymbol g|$ is the relative speed.
Assuming that the particles have unit mass,  
the post-collisional velocities, $\boldsymbol{v}'$ and $\boldsymbol{w}'$,
are given in terms of the pre-collisional velocities, $\boldsymbol{v}$ and 
$\boldsymbol{w}$, by
\begin{equation}
\boldsymbol{v}' = \boldsymbol{v} + \tfrac 12 \, \left( \boldsymbol{g} - g\, \Theta \, \right),
\qquad
\boldsymbol{w}' = \boldsymbol{w} - \tfrac 12\, \left( \boldsymbol{g} - g\, \Theta \, \right),
\end{equation}
for some  direction vector, $\Theta$, on the unit sphere, $S^2$.
We assume that the collisions are modeled using an interparticle potential
of the form, $\phi(r) = r^{-(s-1)}$, for some $1<s\leq \infty$. In this case, 
the collision kernel, $B$, is of the form 
$B(g,\chi) = g^\lambda \widetilde B(\cos \chi)$, where $\lambda=(s-5)/(s-1)$,
and the scattering angle, $\chi$, is given by $\cos \chi = (\boldsymbol{g} \cdot \Theta) /g$. 
 For the main results in this paper, we further 
assume that $0\leq \lambda \leq 1$,
and that the collisions are isotropic, so that $\widetilde B$ is constant. 
The cases $\lambda = 0$ and $\lambda=1$ are those of  Maxwell and hard-sphere collisions, respectively.

Rather than computing  $Q$ itself, Gamba and Tharkabhushanam
consider a truncation, $Q^{\text{tr}}$,
of the collision operator defined by
\begin{equation} \label{eq:Qtr}
Q^{\text{tr}}(\boldsymbol{v})
=\int\limits_{|\boldsymbol g | \leq g_{\text{tr}}} \int\limits_{S^2}
\left[ f(\boldsymbol{v}')f(\boldsymbol{w}')-f(\boldsymbol{v})f(\boldsymbol{v} +\boldsymbol{g}) \right] B\left(g,\frac{\boldsymbol{g} \cdot \Theta}{g}\right) d\Theta d\boldsymbol{g},
\end{equation}
for some choice of truncation parameter, $g_{\text{tr}}$. 
To explain why  it is necessary to
truncate the $\boldsymbol g$-integral in~\eqref{eq:Qtr}, we 
briefly review the derivation of the method.

We define the Fourier transform of a 
function, $F$, on velocity space to be
\begin{equation}
\widehat{F}(\boldsymbol{\zeta}) :=
(2\pi)^{-3/2}\int_{\mathbb{R}^3} 
F(\boldsymbol{v})e^{-i\boldsymbol{\zeta} \cdot \boldsymbol{v}}d\boldsymbol{v}.
\label{QhatDef}
\end{equation}
Using the weak form of the collision operator, 
Gamba and Tharkabhushanam first show that 
\begin{equation}
\widehat{Q}^{\text{tr}}(\boldsymbol{\zeta}) =
(2\pi)^{-3/2} \int_{|\boldsymbol g | \leq g_{\text{tr}}} G(\boldsymbol{g},\boldsymbol{\zeta})
\int_{\mathbb{R}^3} f(\boldsymbol{v})f(\boldsymbol{v} - \boldsymbol{g}) 
\, e^{-i\boldsymbol{\zeta} \cdot \boldsymbol{v}} \, d\boldsymbol{v} d\boldsymbol{g}, 
 \label{fvfvg}
 \end{equation}
 where
\begin{equation}
G(\boldsymbol{g},\boldsymbol{\zeta}) = e^{\frac i2\boldsymbol{\zeta} \cdot \boldsymbol{g}} \int_{S^2}  B\left(g, \frac{\boldsymbol{g}\cdot \Theta}{g}\right) 
e^{-\frac i2 g \,\boldsymbol{\zeta} \cdot \Theta} d\Theta \,\,-\,\, \int_{S^2} B\left(g, \frac{\boldsymbol{g}\cdot \Theta}{g}\right) \, d\Theta.
\label{Formular4G}
\end{equation} 
In the special case of isotropic, inter-particle collisions, the kernel, $G$, is given by
\begin{equation}\label{S2}
G(\boldsymbol{g},\boldsymbol{\zeta})
= 4\pi \widetilde{B} g^{\lambda} \,\left[ e^{\frac i2 \boldsymbol{\zeta} \cdot \boldsymbol{g}}\,\, \text{sinc}\left({g \zeta}/{2}\right) -1\right],
\end{equation}
where $\zeta = | \boldsymbol{\zeta} |$.

Although $G$ is not an integrable function of $\boldsymbol g$ on $\mathbb R^3$,
it does define a tempered distribution~\cite{DuistermaatKolk:2010}. Therefore, 
when $g_{\text{tr}}=\infty$,
the integral~\eqref{fvfvg} for $\widehat Q$  
converges since the velocity pdfs decay exponentially.  However, 
the final step in the derivation of the method 
involves taking the Fourier transform of $G$ with
respect to $\boldsymbol g$, which is a divergent improper integral over
$\mathbb R^3$ when $g_{\text{tr}}=\infty$. This is why it is necessary to truncate
the $\mathbf g$-integral in \eqref{eq:Qtr}. 
Specifically, we define the convolution weighting function by
\begin{equation}\label{GhatDefFinite}
\widehat{G}^\text{tr}(\boldsymbol{\xi},\boldsymbol{\zeta}) = \int_{|\boldsymbol{g}|\leq g_{\text{tr}}} G(\boldsymbol{g},\boldsymbol{\zeta})e^{-i\, \boldsymbol{\xi}\cdot \boldsymbol{g}}\, d\boldsymbol{g}.
\end{equation}
Then, by the convolution theorem, the Fourier transform of the 
truncated collision operator is given by a generalized convolution
integral of the form
\begin{equation} \label{Qhat}
\widehat{Q}^{\text{tr}}(\boldsymbol{\zeta}) = (2\pi)^{-3/2}
\int_{\mathbb{R}^3} \widehat{f}(\boldsymbol{\zeta} - \boldsymbol{\xi}) 
\widehat{f}(\boldsymbol{\xi}) \widehat{G}^{\text{tr}}(\boldsymbol{\xi},\boldsymbol{\zeta}) \, d\boldsymbol{\xi}.
\end{equation}
Since the convolution weighting function is independent of the 
velocity pdfs, $\widehat{G}^{\text{tr}}$ can be precomputed. Therefore
the computational cost of computing the collision operator using \eqref{Qhat} is $\mathcal O(N^6)$.

Haack et al.~\cite{Haack:2013} 
showed that in the special case of isotropic Maxwell collisions, for which
$B=(4\pi)^{-1}$, we have
\begin{equation}
\widehat{G}^\text{tr}(\boldsymbol{\xi},\boldsymbol{\zeta}) \,\,=\,\, 
\widehat{G}_1(\zeta/2,|\boldsymbol{\xi}-\boldsymbol{\zeta}/2|)
-\widehat{G}_2(\xi),
\label{eq:GtrFormula}
\end{equation}
with
\begin{align}
\widehat{G}_1(X,Y)&=\frac{2\pi}{pqXY}\left[q\,\sin(g_{\text{tr}}\,p)-p\,\sin(g_{\text{tr}}\,q)\right],\label{Max:G1}\\
\widehat{G}_2(Z)&= \frac{4\pi }{Z^3}\left[
\sin(g_{\text{tr}}Z)-g_{\text{tr}}\,Z\cos(g_{\text{tr}}Z)\right],
\label{Max:G2}
\end{align}
where
$p=X-Y$ and $q=X+Y$.

\section{Preliminary numerical study}\label{Sec:Prelim}

In this section, we  demonstrate that the choice of the truncation 
parameter, $g_{\text{tr}}$, in~\eqref{eq:Qtr} plays a
critical role in the accuracy of the numerical computation of $Q$.
Clearly, if $g_{\text{tr}}$  is too small
then $Q^{\text{tr}}$ will not be  a good approximation to $Q$. 
On the other hand, 
if $g_{\text{tr}}$ is too large then accurate numerical computation of the 
convolution integral~\eqref{Qhat} is not possible 
since the convolution weighting function, $\widehat{G}^{\text{tr}}$, 
in  \eqref{eq:GtrFormula} is 
a slowly decaying oscillatory function of $\boldsymbol \xi$
whose degree of oscillation increases as $g_{\text{tr}}$ increases.

Gamba and Tharkabhushanam show that if the velocity pdfs are compactly supported,
then there is a value $g_{\text{tr}} < \infty$ so that $Q=Q^{\text{tr}}$. Specifically, they show that if the velocity pdfs are zero outside a box 
$[-L,L]^3 \subset \mathbb{R}^3$, then
$f(\boldsymbol{v})f(\boldsymbol{w}-\boldsymbol{g})=0$
whenever $\boldsymbol{v} \in [-L,L]^3$ and $|\boldsymbol{g}| > g_{\text{max}}$, 
where $g_{\text{max}}=2\sqrt{3} L$. 
Therefore, if we choose $g_{\text{tr}} = g_{\text{max}}$, then
 $Q=Q^{\text{tr}}$, by~\eqref{fvfvg}.
In an unpublished article, Haack~\cite{Haack:2013} instead uses 
$g_{\text{tr}}= L$.
However he provides no explanation for the smaller choice of $g_{\text{tr}}$.
As we now show, the choice of $g_{\text{tr}}$ plays a
critical role in the accuracy of the numerical computation. 
To do so, we compute the collision operator for 
the spherically symmetric, analytical  solution of~\eqref{eq:HBE} 
derived by Bobylev, Krook and Wu~\cite{bobylev1975exact, Max:1976},
which is given by
\begin{equation}
f_{\text{BKW}}(\boldsymbol{v},t) = \frac{e^{-v^2/(2KT)}}{2(2\pi K T)^{3/2}} \left( \frac{5K-3}{K} +\frac{1-K}{K^2}\frac{v^2}{T}\right),
\label{eq:BKW}
\end{equation}
where $v=|\boldsymbol{v}|$ and $K=1-e^{-t/6}$. 
The parameter, $T$, is the temperature, which we set to $T=1$. We compute
the collision operator at the initial time of $t_0=5.5$, which is chosen 
to ensure that $f_{\text{BKW}}>0$. 
Following Haack~\cite{Haack:2013}, we choose the half-width
of the computational domain to be $L=2R$, where $R$ is a measure of the 
effective support for the velocity pdf. Based on a suggestion of
Bobylev and Rjasnow~\cite{Bobylev:1997}, Haack chooses $R=2\sqrt{2}T$, where
$T$ is the temperature of the distribution, which results in $L\approx 5.66$.
This choice is justified
by the observation that the velocity 
pdf typically decreases as does $\exp(-v^2/2T)$ for large $v$.

\begin{figure}[t!]
\centering
\includegraphics[width=0.49\textwidth]{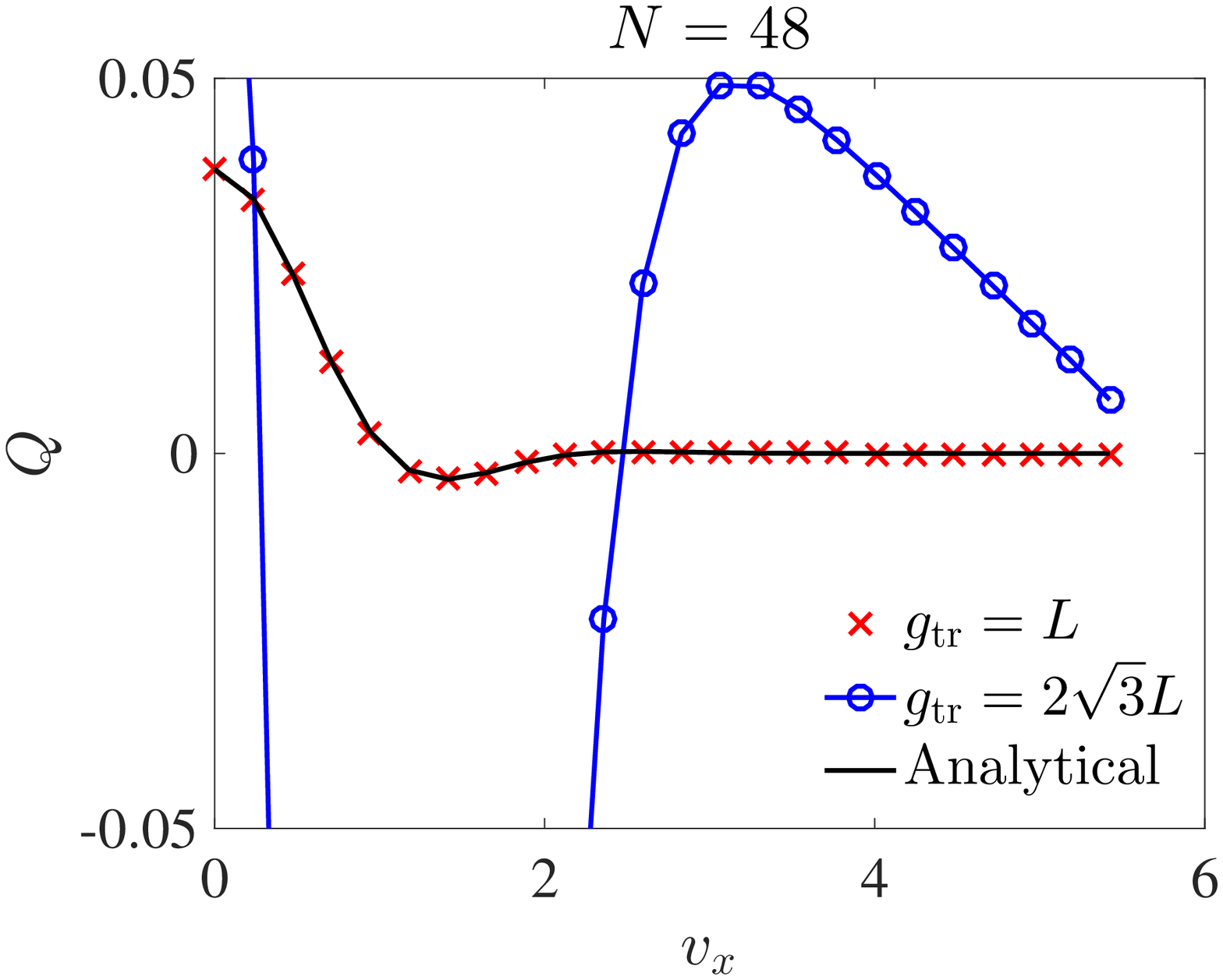}
\includegraphics[width=0.49\textwidth]{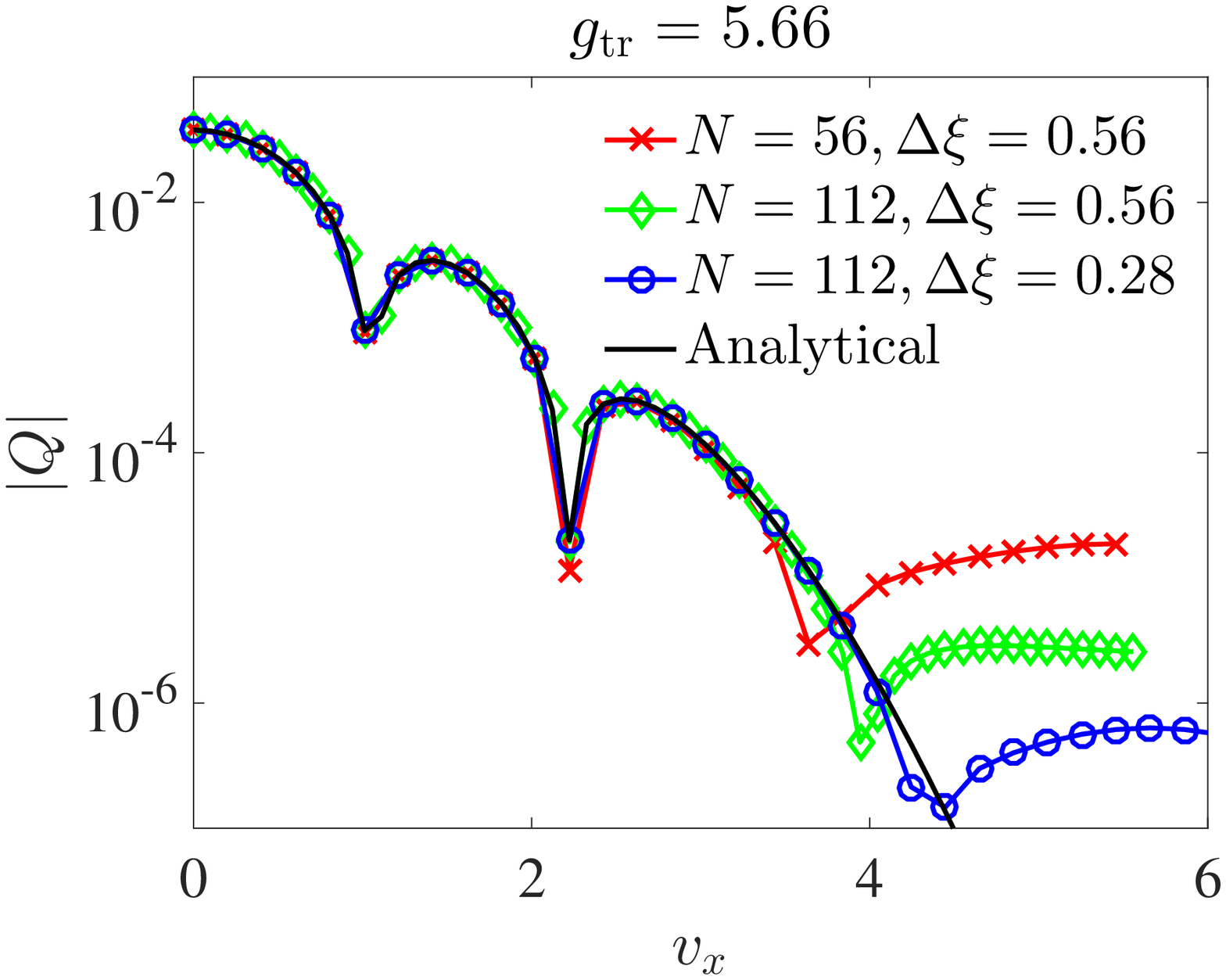}\\
\caption{Slices in the $v_x$-direction of the collision operator at $(v_y,v_z)=(0,0)$
for the BKW pdf \eqref{eq:BKW} with $T=1$ at $t=5.5$.
Left: Linear scale plots of $Q$ for the values of 
$g_{\rm tr}$ obtained using the formulae given in \cite{GambaThark228p2012} and 
\cite{Haack:2013}
with $N=48$ discretization points in each velocity direction.
 Right: Log scale plots of $|Q|$ 
for the value of $g_{\rm tr}$ given in \cite{Haack:2013} for
different values of  $N$ and $\Delta\xi$.  In both panels, 
the black solid curves show the results obtained using the analytical formula~\eqref{eq:BKW}.
}
 \label{fig:BadQ}
\end{figure}

In Fig.~\ref{fig:BadQ} (left) we plot slices of the collision operator
on a linear scale  with $g_{\text{tr}}= L$ and $g_{\text{tr}}= 2\sqrt{3}L$.
For these results we used $N=48$ discretization points in each
velocity direction. %%The results with $N=24$ were essentially the same.
The result with $g_{\text{tr}}= L$ agrees well with the analytical
formula for $Q$ obtained using \eqref{eq:HBE} and \eqref{eq:BKW}.
However, the result with $g_{\text{tr}}= 2\sqrt{3}L$ is 
far from being correct, which suggests that in their numerical simulations
Gamba and Tharkabhushanam actually used a significantly smaller
value of $g_{\text{tr}}$, although they do not state which value they chose. 

The reason for the lack of agreement with $g_{\text{tr}}= 2\sqrt{3}L$ is that
the integral  \eqref{Qhat} for $\widehat Q^{\text{tr}}$ 
is numerically computed 
using values of $\widehat f$ and $\widehat G^{\text{tr}}$ 
on a grid with spacing 
\begin{equation}%\label{deltaXi}
\Delta \xi = \frac{\pi}{L}=\frac{2 \sqrt{3}\pi}{g_{\text{tr}}}.
\end{equation}
This grid spacing is determined by the standard discretization of the
Fourier transform.
With this grid spacing there can be significant  error in the 
numerical computation of the integral, since as we see from \eqref{Max:G1} and
\eqref{Max:G2},
the kernel, $\widehat G^{\text{tr}}$,
oscillates on a length scale of approximately $2 \pi/g_{\text{tr}}$.
Increasing the number of discretization points, $N$, in each velocity dimension 
does not change $\Delta \xi$. 
On the other hand, the agreement is much better 
with $g_{\text{tr}}= L$, since the 
frequency of oscillation of $\widehat G^{\text{tr}}$ is 
approximately half that of the sampling frequency, $\Delta \xi = \pi/g_{\text{tr}}$,
in accord with the Nyquist-Shannon sampling theorem. 
 
Because we are interested in computing the low-probabilty tails of the
velocity pdf, it is important to determine how the choice of $g_{\text{tr}}$ affects the
relative error in the numerically computed values of collision operator
at large speeds, $v = | \mathbf v|$.
To start investigating this question, 
in Fig.~\ref{fig:BadQ} (right) we plot slices of the absolute value of the collision operator
on a logarithmic scale. 
We note that the  cusps evident in the log-scale plots occur where $Q$ changes sign.
The three numerical results were obtained 
using $g_{\text{tr}}= 5.66$ but with different choices for $N$ and for
the grid spacing $\Delta \xi=\pi/L$ used in the numerical
computation of the integral~\eqref{Qhat}. The red line with crosses and the green line with diamonds shows the results with $N=56$ and $N=112$, respectively. In both cases we choose $L=g_{\text{tr}}=5.66$, which results in $\Delta \xi= 0.56$. 
We see from these results that doubling $N$ increases the accuracy of the computation by about an order of magnitude. The reason is that $\widehat G^{\text{tr}}$
decays slowly as $|\boldsymbol\xi|\to\infty$ and increasing $N$ increases
the size of the domain of integration in frequency space.
The blue curve with circles shows the result with $N=112$ and 
$L=2g_{\text{tr}}=11.32$, so that $\Delta \xi =0.28$.
Comparing the red and blue curves, we see that 
halving $\Delta \xi$ 
increases the accuracy of the computation by about two orders of magnitude.
The reason is the smaller grid spacing better captures
the oscillations of $\widehat G^{\text{tr}}$.

One of our main goals in this paper is to  fix $L$ and investigate how small we can
choose both $g_{\text{tr}}$ and $N$ so as to accurately compute the
velocity pdf down to a desired probability level. 
If we let $Q^{\text{NC}}$ denote the collision operator
obtained by numerically computating of $Q^{\text{tr}}$, 
 then the total error,
 \begin{equation}
 \mathcal E_{\text{tot}} \,\,:=\,\, \left|Q - Q^{\text{NC}}
 \right|,
 \label{TotErr}
 \end{equation}
is bounded by
\begin{equation}\label{errorsplitformula}
 \mathcal E_{\text{tot}}
 \,\,\leq\,\, 
\left|Q-Q^{\text{tr}}\right| \,\,+\,\,
\left|Q^{\text{tr}}-Q^{\text{NC}}\right|.
\end{equation} 
The first term on the right hand side of \eqref{errorsplitformula} is the error inherent in the 
truncation of the collision operator, and the second term is the error
in the numerical computation of the truncated  operator. 
Since in practice, we have limited computational resources
the choice of $g_{\text{tr}}$ involves a trade off between these two sources of error.

\section{An error estimate for the truncated collision operator}\label{Sec:Tr}
In this section, we first review a theorem of 
Cercignani et al.~\cite{cercignani2013mathematical}
on the convergence of  $Q^{\text{tr}}$ to $Q$
as $g_{\text{tr}}\to\infty$. However, since this theorem does not
include an error bound, it is of limited utility for numerical
computation. Then, we derive an upper bound on  the pointwise
error between $Q$ and $Q^{\text{tr}}$,
assuming that both operators are computed exactly.
This bound  provides some additional theoretical justification for the spectral-Lagrangian
method. In Section~\ref{Sec:Results}, we will use this upper bound 
to guide the choice of $g_{\text{tr}}$ in numerical computations.

Cercignani et al.~\cite{cercignani2013mathematical} use 
the truncated collision operator in a proof of  an existence and uniqueness
theorem for solutions of the homogeneous Boltzmann equation. 
In their proof they consider the initial value problem, 
\begin{align}\label{eq:HBEtr}
\begin{aligned}
\frac{\partial f^M}{\partial t}\,\,&=\,\, Q^M(f^M,f^M), \\
f^M(0,\cdot) \,\,&=\,\, f_0,\\
\end{aligned}
\end{align}
where $Q^M = Q^{\text{tr}}$ with $g_\text{tr} = M < \infty$.
They show that \eqref{eq:HBEtr}
has a unique nonnegative solution $ f^M\in C^1( [0,T],L^1(\mathbb R^3))$ 
for all $T > 0$, provided that $f_0\in L^1(\mathbb R^3) $ is nonnegative. 
In addition, they show that the total mass, momentum, and energy are conserved
by \eqref{eq:HBEtr}. 
Applying the Dunford-Pettis theorem to the set $\{f^M\}$ in 
$C^1( [0,T],L^1(\mathbb R^3))$, they extract a 
weakly convergent subsequence $f^n \rightharpoonup f$, with $f\in C^1( [0,T],L^1(\mathbb R^3))$ nonnegative, and prove that
$Q^{M_n} (f^n,f^n) \rightharpoonup Q(f,f)$. 
Here, by weak convergence we mean
convergence of the sequence obtained after integration against a test function
in $L^\infty (\mathbb R^3)$. 

To assess the trade off discussed at the end of Section~\ref{Sec:Prelim}, we now
present an upper bound for the error inherent in the truncation of $Q$,
\emph{i.e.}, for the first term on the right hand side of \eqref{errorsplitformula}. 
Specifically, we let
\begin{equation}\label{eq:trerr}
\mathcal E_{\text{tr}}(\mathbf v) \,\,:=\,\,
\left|Q(\mathbf v)-Q^{\text{tr}}(\mathbf v)\right|,
\end{equation}
where $Q = Q(f,f)$ and 
$Q^{\text{tr}} = Q^{\text{tr}}(f,f)$.
Since these two collision operators are defined in terms of the same velocity pdf, $f$, 
this result has a different
character than that of Cercignani.

To obtain this result, rather than assuming that
the support of $f$ is compact, we instead assume that 
$f$ is  bounded above by a Maxwellian pdf. 
For Maxwell-type collisions between particles of the same type, Bobylev and Gamba~\cite{BobylevGamba:2017} proved 
that, if the initial condition satisfies 
$f_0(\boldsymbol v) \leq c_0 \,e^{- k_0 v^2}$, then there are constants $c\geq c_0$ 
and  $k\leq k_0$ so that $f(\boldsymbol v,t) \leq c \,e^{- k v^2}$ for all $t>0$. Moreover, they
provide formulae for $c$ and $k$ in terms of $c_0$, $k_0$, and the initial pdf $f_0$.  
Gamba et al.~\cite{GambaPanferov:2007} proved similar results for other 
interparticle collision kernels. 
Consequently, the assumptions we make in the
following theorem and in its corollary are reasonable.   

\begin{theorem}
\label{thmerror}
 Suppose that  $f(\boldsymbol{v})\leq c  \,e^{- k v^2}$, where $v=|\boldsymbol{v}|$, and 
that
the collision kernel is of the form $B(g, \chi) = g^\lambda\, \widetilde{B}$, where $\widetilde{B}$ is a positive constant and $0\leq \lambda \leq 1$. Then, the 
error \eqref{eq:trerr} in the truncation  of the collision operator is bounded by
\begin{equation}\label{anarelerror}
\mathcal E_{\operatorname{tr}}(\mathbf v)
 \,\,\leq\,\,
 \mathcal E^{\operatorname{UB}}_{\operatorname{tr}}(g_\text{tr},v)
 := c \exp(- k v^2)\,
 \mathcal E_{\operatorname{rel}}(g_{\operatorname{tr}},{v}),
\end{equation}
where
\begin{equation}
\mathcal E_{\operatorname{rel}}(g_{\operatorname{tr}},{v}) \,\,:=\,\,
16 {\pi}^2  \widetilde{B} \, c\, 
 \int_{g=g_{\operatorname{tr}}}^\infty e^{-{k} (v-g)^2} \left[  \frac{1-e^{-4{k}  v g}}{4 {k} v g}\right] g^{\lambda+2} \, dg.
\label{Error_Estimate}
\end{equation}
\end{theorem}

Since $Q$ is the rate of change of $f$,
we expect $Q(\mathbf v)$ to be on the order
of $c\exp(- k v^2)$ or less.
Therefore, we can regard the error, $\mathcal E_{\operatorname{rel}}$, in~\eqref{Error_Estimate}
as a measure of the relative error between $Q$ 
and $Q^{\text{tr}}$. Note that
since $Q$ has zeros, we have defined this error to be relative to a 
Maxwellian pdf rather than to $Q$.
The upper bound, 
$\mathcal{E}_{\operatorname{rel}}$, in \eqref{Error_Estimate} can be used to 
guide the choice of $g_{\text{tr}}$ in numerical simulations by ensuring that
 $\mathcal E_{\operatorname{rel}}(g_{\operatorname{tr}},{v}) < 10^{-m}$ 
 for a desired value of $m$, over a given range of values for $v$.

\begin{proof}
Using the assumptions in the statement of the theorem,
\begin{equation}
|Q(\boldsymbol v) -Q^{\operatorname{tr}}(\boldsymbol v)|
\,\,\leq\,\,  
\widetilde B \max \left\{ \mathcal E_1 , \mathcal E_2 \right\},
\label{maxerrorinequality}
\end{equation}
where 
\begin{equation}\label{E1}
\mathcal E_1 = \int\limits_{|\textbf{g}| \geq g_{\text{tr}}}  \int_{S^2}
 f\left(\boldsymbol{v} + \tfrac{1}{2} \, \left( \boldsymbol{g} - g\, \Theta \, \right)\right)f\left(\boldsymbol{v}+ \tfrac{1}{2} \, \left( \boldsymbol{g} + g\, \Theta \, \right) \right) g^\lambda \, d\Theta\, d\boldsymbol{g},
\end{equation}
and 
\begin{equation}\label{E2}
\mathcal E_2 =  \int \limits_{|\boldsymbol{g}|\geq g_{\text{tr}}} \int_{S^2}
f(\boldsymbol{v})f(\boldsymbol{v}+\boldsymbol{g}) \,g^\lambda 
\, d\Theta\, d\boldsymbol{g}.
\end{equation} 
The inequality \eqref{maxerrorinequality} holds since $\mathcal E_1$ and $\mathcal E_2$ are both positive. Using the upper bound we have assumed for $f$, we find that
\begin{align}
\mathcal E_2 
&\,\,\leq \,\, 
4 \pi c^2 e^{-2k v^2} \int \limits_{|\boldsymbol{g}|\geq g_{\text{tr}}}  e^{-k(2\boldsymbol{v} \cdot \boldsymbol{g}+g^2)} \,g^\lambda \, d\boldsymbol{g} 
\label{E2A}
\\
&\,\,=\,\,
8 \pi^2 c^2 e^{-2k v^2} 
\int_{g_{\text{tr}}}^{\infty}  e^{-k g^2} g^{\lambda +2}\, \int_0^{\pi}  
e^{-2\,k v g\cos \phi} \, \sin \phi \,d\phi\,  dg \nonumber
\\
&\,\,=\,\, 16 {\pi}^2 \, c^2 \,e^{-k v^2} \int_{g_{\text{tr}}}^\infty e^{-k(v-g)^2} \left[  \frac{1-e^{-4k v g }}{4 k v g}\right] g^{\lambda+2} \, dg.
\label{E2B}
\end{align}

Similarly,  
\begin{align}
\mathcal E_1 &\,\,\leq \,\,
c^2
\int \limits_{|\boldsymbol{g}|\geq g_{\text{tr}}} \int_{S^2}  e^{-k |\boldsymbol{v} + \frac{1}{2}(\boldsymbol{g}-g \Theta)|^2} \,  e^{-k|\boldsymbol{v} + \frac{1}{2}(\boldsymbol{g}+g \Theta)|^2} \, g^\lambda d \Theta \, d \boldsymbol{g}
\\
&\,\,\leq\,\,
 4\pi c^2 \, e^{-2{k} v^2}  \int \limits_{|\boldsymbol{g}|\geq g_{\text{tr}}} e^{-{k} (2\boldsymbol{v} \cdot \boldsymbol{g}+g^2)} g^\lambda 
  \, d \boldsymbol{g}.
\end{align}
As in \eqref{E2A} and \eqref{E2B}, we find that
\begin{equation}\label{E1estimatefinal}
\mathcal E_1 \leq 16 {\pi}^2 \, c^2 \,e^{-{k} v^2} \int_{g_{\text{tr}}}^\infty e^{-{k}(v-g)^2} \left[  \frac{1-e^{-4{k}  v g}}{4 {k} v g}\right] g^{\lambda+2} \, dg .
\end{equation} 
The required estimate now follows from \eqref{maxerrorinequality}, \eqref{E1estimatefinal},
and the fact that the right hand side of \eqref{E2B} is bounded above 
by the  right hand side of \eqref{E1estimatefinal}.
\end{proof}

\begin{corollary}\label{corerror}
Suppose that the  assumptions of Theorem~\ref{thmerror} hold and that $\lambda=0$.
Then, provided $k$ and  $g_{\operatorname{tr}}$ are both large enough,
we have the asymptotic formulae
\begin{equation}\label{eq:Asymp3}
\mathcal E_{\operatorname{rel}} (g_{\operatorname{tr}},v) 
\,\,\approx \,\,
\begin{cases}
c (\frac{\pi}{k})^{3/2} & \qquad \text{if } g_{\operatorname{tr}} < v, \\

\frac{1}{2} [ c (\frac{\pi}{k})^{3/2} \,\,+\,\, \frac{1}{kg_{\operatorname{tr}}} ]
& \qquad \text{if } g_{\operatorname{tr}} = v, \\

\frac{\pi c}{2k^2} \, \frac{ e^{-k(g_{\operatorname{tr}} - v)^2}}{g_{\operatorname{tr}}-v}
\, \frac{g_{\operatorname{tr}}}{v} & \qquad \text{if } g_{\operatorname{tr}} > v.

\end{cases}
\end{equation}
In particular, if $g_{\operatorname{tr}} \gtrsim v \gg 0$, then
\begin{equation}\label{eq:Asymp1}
\mathcal E_{\operatorname{rel}} (g_{\operatorname{tr}},v) 
\,\,\approx\,\, 
\frac{\pi c}{2k^2} \, \frac{ e^{-k(g_{\operatorname{tr}} - v)^2}}{g_{\operatorname{tr}}-v},
\end{equation}
which is a rapidly decaying function of $g_{\operatorname{tr}}-v$. 
\end{corollary}

Equation~\eqref{eq:Asymp1} 
yields the following simple strategy for choosing the parameter, $g_{\operatorname{tr}}$, 
in numerical computations of the low-probability tails. Specifically, 
to guarantee $Q^{\operatorname{tr}}$
is an accurate approximation to $Q$ at $\mathbf v$, we should choose
$g_{\operatorname{tr}}$ to be slightly larger than $v$. 
As we will see in Section~\ref{Sec:Results}, exactly how much larger depends on the values of $c$ and $k$ and the desired degree of accuracy. 
%This result also provides a theoretical justification for Haack's choice of 
%$g_{\operatorname{tr}}=L$~\cite{Haack:2013}. 
%However,  sometimes we can use a smaller value of $g_{\operatorname{tr}}$
%since the upper bound, $\mathcal E^{\operatorname{UB}}_{\operatorname{tr}}$,
%in Theorem~\ref{thmerror} may be suboptimal. 

\begin{proof}
Under the assumptions of \eqref{corerror}, 
\begin{equation}\label{SimpleE}
\mathcal E(g_{\operatorname{tr}},v) \,\,=\,\, 4\pi c
\int\limits_{g_{\operatorname{tr}}}^\infty 
e^{-k(v-g)^2}\, \frac{1-e^{-4kvg}}{4kvg}\, g^2\, dg.
\end{equation}
In the cases that $g_{\operatorname{tr}} \leq v$, 
we apply Laplace's method~\cite{bleistein1986asymptotic} as follows.
First, recall that if a function $\phi:[a,b]\to\mathbb R$ has a single critical point at 
an interior point, $t_0 \in (a,b)$, which is the  absolute minimum of $\phi$,  then
for any sufficiently smooth function, $f$, 
\begin{equation}
\int_a^b e^{-k \phi(t)}\, f(t)\, dt
\,\,\approx\,\, 
\sqrt{\frac{2\pi}{k\phi''(t_0)}} \,\, e^{-k  \phi(t_0)} \, f(t_0), 
\qquad \text{as } k\to\infty.
\end{equation}
When  $g_{\operatorname{tr}} < v$, the result follows by setting 
$t_0=v$, and using the estimate $f(v) = \frac{c\pi}{k} (1 - e^{-4kv^2}) 
\approx \frac{c\pi}{k}$, provided $g_{\operatorname{tr}}$ is large enough. 
When $g_{\operatorname{tr}} = v$, the critical point,
$t_0 = g_{\operatorname{tr}}$, of $\phi$ is an endpoint of the interval
of integration. Then by \cite[(5.1.17)]{bleistein1986asymptotic},
we find that
\begin{equation}
\mathcal E(g_{\operatorname{tr}},g_{\operatorname{tr}})
\,\,\approx\,\,
\frac{\pi c}{k g_{\operatorname{tr}}} \left[
\sqrt{\frac{\pi}{4k}}\,g_{\operatorname{tr}}(1- e^{-4kg_{\operatorname{tr}}^2})
+ \tfrac 1{2k} 
\left( 1 - (1+4kg_{\operatorname{tr}}^2) e^{-4kg_{\operatorname{tr}}^2}\right)
\right].
\end{equation}
The result now follows provided $g_{\operatorname{tr}}$ is large enough. 

Finally, in the case  $g_{\operatorname{tr}}>v$, the change of variables
$t=(g-v)^2 - (g_{\operatorname{tr}}-v)^2$ transforms \eqref{SimpleE} to 
\begin{equation}
\mathcal E(g_{\operatorname{tr}},v) \,\,=\,\, 2\pi c \,e^{-k(g_{\operatorname{tr}}-v)^2}
\, \int_0^\infty e^{-kt} F(\sqrt{t+ (g_{\operatorname{tr}}-v)^2}\,;\,k,v)\, dt,
\end{equation}
where
\begin{equation}\label{F}
F(u\,;\,k,v) = \frac{1-e^{-4kv(v+u)}}{4kvu}\,\,\frac{(v+u)^2}{u}.
\end{equation}
The result now follows from Watson's Lemma~\cite{bleistein1986asymptotic},
which states that in the limit as $k\to\infty$, we have that
 $\int_0^\infty e^{-kt}H(t)\, dt $ $\approx H(0)/k$. Although~\eqref{eq:Asymp3} is only guaranteed to hold in the limit $k\to\infty$,
in Section~\ref{Sec:Results} we will show it that is quite accurate even for 
$k=\mathcal{O}(1)$.
\end{proof}

\section{Numerical Method}\label{Sec:Method}
We implement Gamba's method 
as in~\cite{GambaThark228p2012,Haack:2013}. 
The computational grids in velocity and Fourier space are
 defined in terms of a maximum speed, $L$, and
the number of grid points, $N$, in each dimension.
Unless otherwise noted, we use $L=10$. 
We represent the velocity pdf on a domain $[-L,L]^3\in\mathbb R^3$
using a regular grid, $v_k=-L+k\Delta v$,
in each velocity dimension, where $\Delta v = 2L/N$. The corresponding
domain in Fourier space is $[-\zeta_{\rm max},\zeta_{\rm max}]^3$, where
$\zeta_{\rm max} = N\pi/2L$, with grid points,
$\zeta_m = -\zeta_{\rm max} + m \Delta \zeta$, where $\Delta \zeta = \pi/L$.
We calculate the integral \eqref{Qhat} using the
trapezoid rule, which---like the discrete
Fourier transform---is spectrally accurate for functions that decay rapidly
at the boundary of the computation domain~\cite{trefethen2000spectral}. 
In addition, for the  numerical results in 
Sections~\ref{Sec:Maxwellian}--\ref{Sec:MixMax} below we enforce
conservation of the density, momentum, and energy using the
Lagrangian projection method of~\cite{GambaThark228p2012}, which
 amounts to projecting the collision operator
onto a linear subspace in $L^2(\mathbb R^3)$. 
However, we found that the results
in these subsections are visually indistinguishable from 
those we obtained without the application of the 
Lagrangian projection method.
Because of the large cost of computing the collision operator, 
we use a multistep method to solve the 
system of differential equations corresponding to \eqref{eq:HBE}, which 
allows us to take larger time steps resulting in a  fewer total number of function evaluations. 
Specifically, we used
the fourth order Adams-Bashforth method~\cite{Kincaid:1991}
\begin{equation}\label{Adams}
f_{i+4} = f_{i+3} + \frac{\Delta t}{24} \Big[55 Q_{i+3}-59 Q_{i+2}+37 Q_{i+1} -9 Q_i\Big],
\end{equation}
where $f_i(\boldsymbol v) = f(t_i,\boldsymbol v)$ and 
$Q_i = Q^{\operatorname{tr}}(f_i,f_i)$.
To initialize this multistep method we used the 
fourth order Runge Kutta method to compute the solutions at the first four time steps.  

\section{Numerical Results}\label{Sec:Results}
In this section we present the results of the numerical simulations we performed to
test the limits of the spectral-Lagrangian method. We show results
for the following choices of initial condition: a Maxwellian, the spherically symmetric, analytical  solution of  Bobylev, Krook and Wu~\cite{bobylev1975exact, Max:1976}, a cylindrically symmetric
modification of the BKW initial condition, and two mixtures of Maxwellians.
We study the convergence of the numerically computed truncated collision operator,
$Q^{\operatorname{NC}}$, to the collision  operator, $Q$, 
validate the bound we obtained for the relative error between 
$Q^{\operatorname{tr}}$ and $Q$, and 
compute the evolution of the velocity pdf to the equilibrium 
Maxwellian distribution. In these simulations our focus is on 
the accuracy with which the velocity pdfs can be computed in the low probability tails. 
Finally, we employ  a simple model of a
plasma  to study
the evolution of the velocity pdf of the electrons 
under the influence of an electron gun source,
electron-electron collisions, and loss into a boundary sheath layer.

The simulations were performed on a 2.0~GHz Intel Xeon processor with
2 CPU's and 14 cores per CPU. The total simulation time (number of cores
$\times$ time per core) for a single computation of the collision operator
ranged from 24~seconds for $N=24$ to 4.5~hours for $N=72$, and scaled
according to the theoretical $\mathcal O(N^6)$ cost.

\subsection{The Maxwellian solution}\label{Sec:Maxwellian}

If the initial velocity pdf is a Maxwellian, 
$f(\mathbf v) = (2\pi T)^{-3/2} \exp(-v^2/2T)$, then the collision operator is identically
zero, $Q\equiv 0$. In Table~\ref{table:MaxMaxwellianQBeforeLP}, with $T=1$, we
plot the $L^\infty$-error in the numerically computed truncation operator, $Q^{NC}$,
for several pairs of values of $N$ and $g_{\operatorname{tr}}$. 
When $g_{\operatorname{tr}}= 4$ and 8, the error decreases to the
level of the round-off error for the Fourier transform as $N$ increases from 24 to 72.
However, when $g_{\operatorname{tr}}= 12$, the error does not  converge to zero
since the convolution weighting function, ${\widehat G}^{\operatorname{tr}}$,
oscillates on a length scale that is close to $\Delta \xi=\Delta\zeta$. 
The results are significantly worse
when $g_{\operatorname{tr}}$ is increased to 16 and 20. 
These results show that, if $N$ is large enough 
to capture the  slow decay of  ${\widehat G}^{\operatorname{tr}}$ ($N\geq 48$), 
then we can obtain a large gain in 
the accuracy of the generalized convolution integral~\eqref{Qhat} for 
${\widehat Q}^{\operatorname{tr}}$ by 
choosing $g_{\operatorname{tr}}<L$, thereby reducing the
oscillation of ${\widehat G}^{\operatorname{tr}}$ relative to the
grid spacing in Fourier space. 
For more general initial conditions, because $\widehat f$ is typically
smoother than ${\widehat G}^{\operatorname{tr}}$, we expect a similar
rate of convergence of ${Q}^{\operatorname{tr}}$ to
${Q}^{\operatorname{NC}}$, i.e., for the second term in \eqref{errorsplitformula}. 
In the next subsections, for several choices of initial condition, we use the error bound 
in Theorem~\ref{thmerror}
 to determine values of $g_{\operatorname{tr}}$ for which we can
guarantee that the error in the second term in \eqref{errorsplitformula} is 
below a given threshold out to a given value of $v$. 

\renewcommand{\arraystretch}{1.2}
\begin{table}[!h]
\centerline{ 
\begin{tabular}{|c||c|c|c|c|c|}
\hline
&$g_\text{tr}=4$ & $g_\text{tr}=8$ & $g_\text{tr}=12$ & $g_\text{tr}=16$ & $g_\text{tr}=20$ \\
\hline
\hline
$N=24$ & $2\times 10^{-5}$      & $3\times 10^{-5}$ &   $4\times 10^{-5}$       & $2\times 10^{-4}$ &   $2\times10^{-1}$ \\ \hline
$N=36$ & $2\times 10^{-9}$      & $4\times 10^{-9}$ &   $4\times 10^{-9}$       & $2\times 10^{-4}$ &   $2\times 10^{-1}$ \\ \hline
$N=48$ & $8\times 10^{-15}$     & $1\times 10^{-14}$ & $5\times 10^{-10}$       & $2\times 10^{-4}$ &   $2\times 10^{-1}$ \\ \hline
$N=72$ & $5\times 10^{-17}$     & $5\times 10^{-17}$ & $5\times 10^{-10}$       & $2\times 10^{-4}$ &   $2\times 10^{-1}$ \\
\hline
\end{tabular}
}  \vskip10pt
\caption{$L^\infty$-error in $Q^{NC}$ for the Maxwellian velocity pdf 
 for different values of $N$ and $g_{\operatorname{tr}}$.}
 \label{table:MaxMaxwellianQBeforeLP}
\end{table}

\subsection{The BKW solution}\label{Sec:BKW}

In this subsection, we compare the results obtained using the numerical method
to the  analytical  solution, $f_{\text{BKW}}$, of  Bobylev, Krook and Wu~\cite{bobylev1975exact, Max:1976} given in \eqref{eq:BKW} with $T=1$. 

We begin by using Theorem~\ref{thmerror} to 
select an appropriate value of the truncation parameter, $g_{\rm tr}$,
for the  velocity pdf, $f_{\text{BKW}}$, at the initial time of $t=5.5$.
We consider two methods for selecting the Maxwellian upper bound
required to apply the theorem. For Method~I we 
choose the width parameter, $k$, in Theorem~\ref{thmerror} to agree with the 
width of the Maxwellian pdf to which the initial condition
converges as $t\to\infty$. This method gives $k = 3/2E$, where
$E$ is the (initial) energy. We then choose the parameter
$c$ to ensure that the resulting Maxwellian pdf is an upper bound for
the velocity pdfs at the initial time. 
For the BKW pdf, Method~I gives $k=0.5$ and $c=0.1$. 
For Method~II, we use the tightest upper Maxwellian bound we could
find for the initial velocity pdf, which resulted in $k=0.8$ and $c=1$.

\begin{figure}[t!]
\centering
\includegraphics[width=0.32\textwidth]{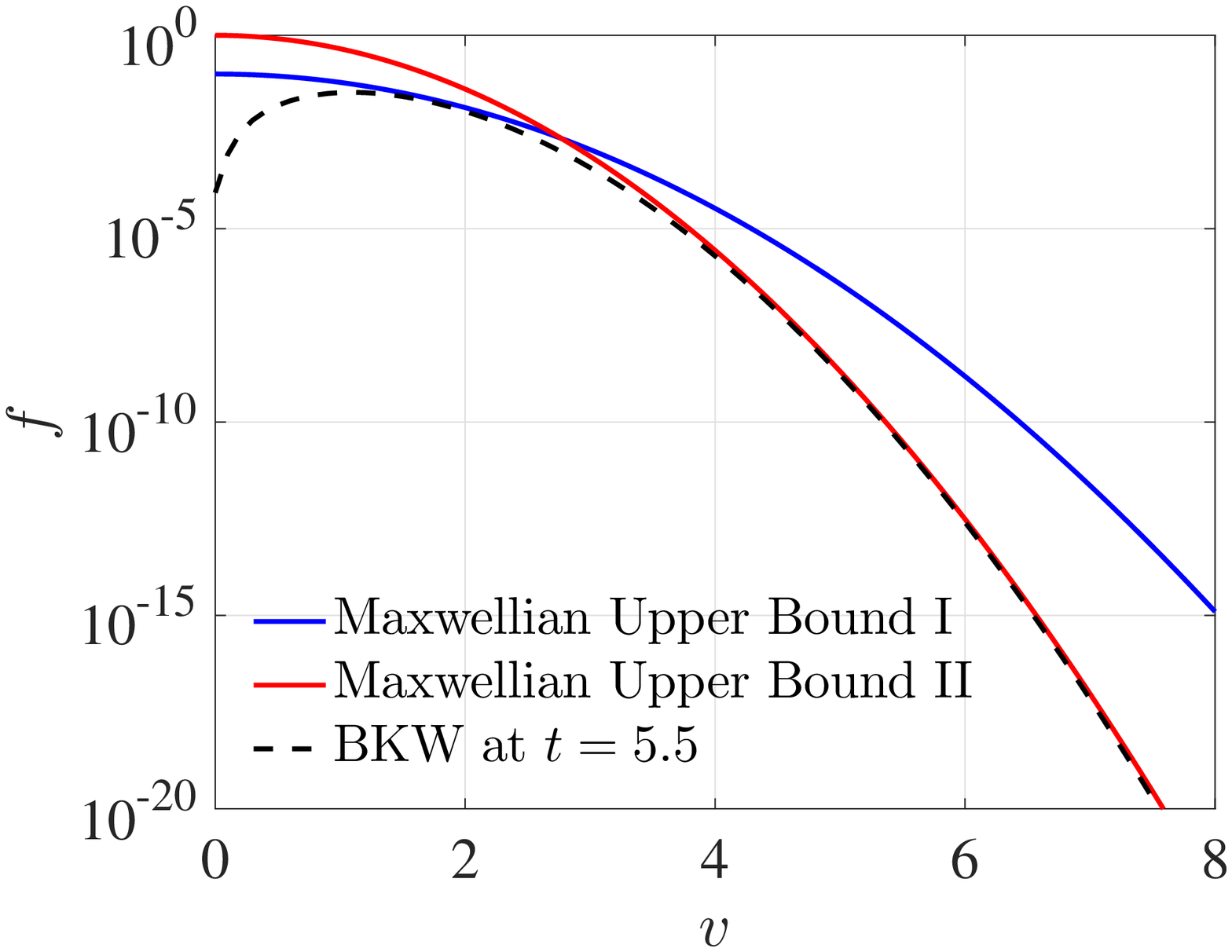}
\includegraphics[width=0.32\textwidth]{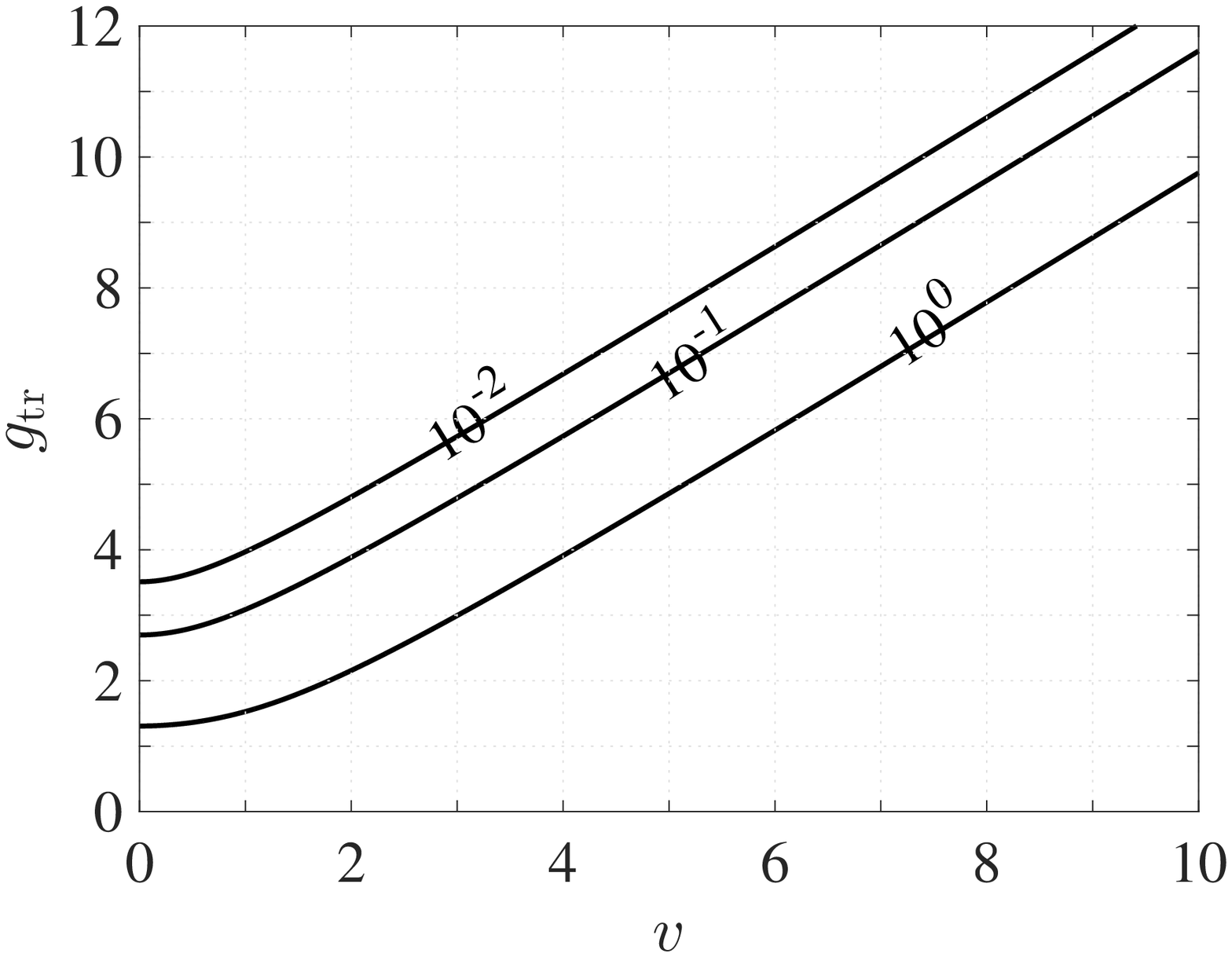} 
\includegraphics[width=0.32\textwidth]{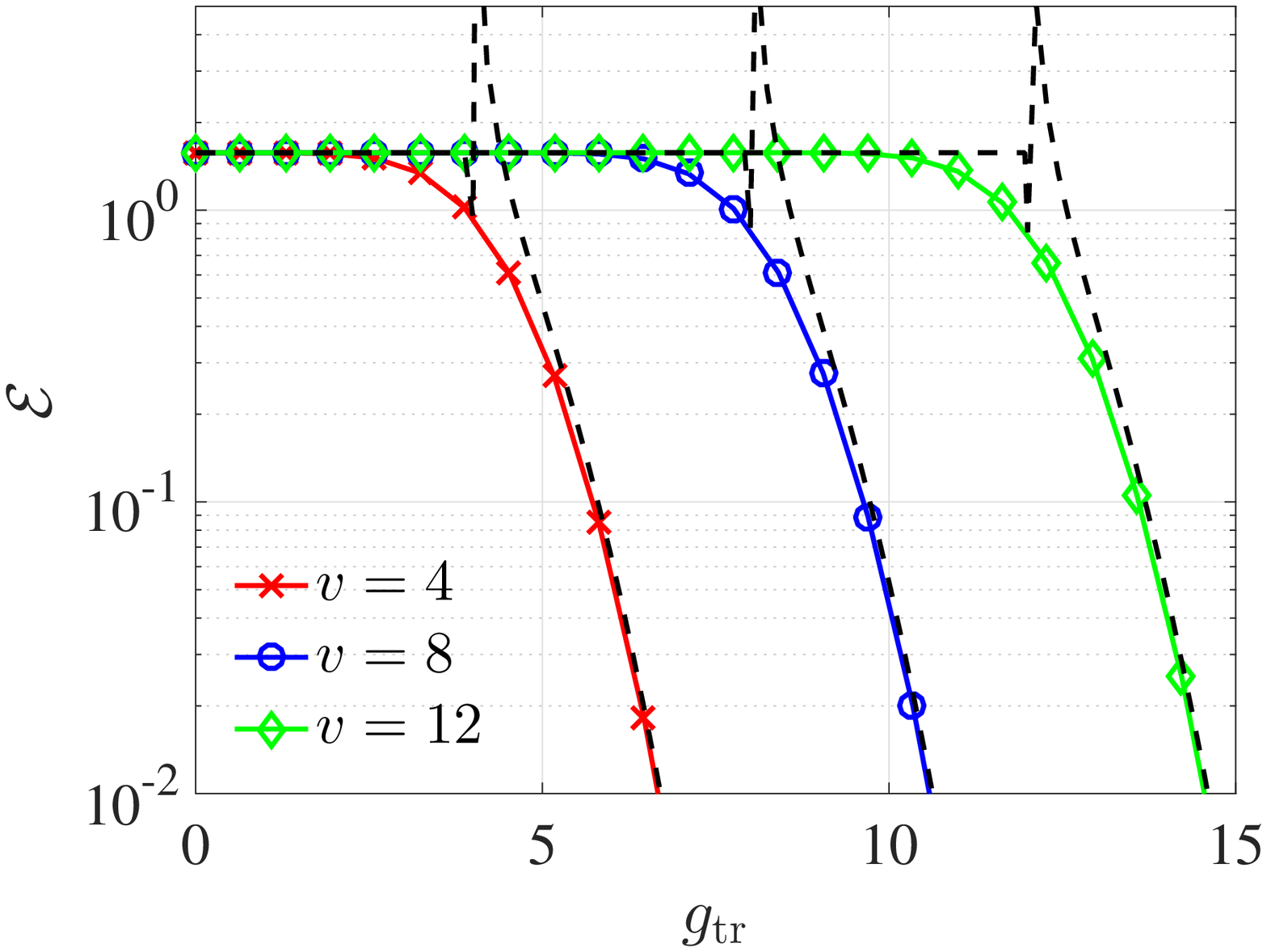}
\caption{
Left: Log-scale plot of the initial BKW velocity pdf given by~\eqref{eq:BKW}
at the initial time, $t=5.5$,
(dashed black curve), together with the
Maxwellian upper bounds obtained using Method~I (solid blue curve) and 
Method~II (solid red curve). 
Middle: Contour plot of the upper bound, $\mathcal E_{\operatorname{rel}}$, for the relative error in the truncation of the collision operator given by \eqref{Error_Estimate},
as a function of speed, $v$, and truncation parameter, $g_{\rm tr}$.
This result was obtained using the Maxwellian upper bound obtained using Method~I.
Right: Slices of 
$\mathcal E_{\operatorname{rel}}(v,g_{\text{tr}})$ for three values of $v$.
The colored solid curves with symbols show the results obtained with
 \eqref{Error_Estimate} and the black dashed curves show the corresponding 
 asymptotic formulae in \eqref{eq:Asymp3}.
}
 \label{fig:RelErrorContours}
\end{figure}

In Fig.~\ref{fig:RelErrorContours} (left), we plot the initial BKW pdf
and the two Maxwellian bounds, and in Fig.~\ref{fig:RelErrorContours} (middle)
we show a contour plot of the bound, $\mathcal E_{\operatorname{rel}}$, 
for the  relative error in the truncation of the collision operator 
given by~\eqref{Error_Estimate}, as a function of  $v$ and $g_{\text{tr}}$. 
For this contour plot we have used the Maxwellian upper bound given by Method~I.
The results obtained with Method~II are quite similar: For each $v$ the 
contours are shifted up or down by about 0.5 in $g_{\text{tr}}$.
The contour plot shows that if we choose $g_{\text{tr}}=6$ then
$\mathcal E_{\operatorname{rel}} < 10^{-1}$ for $v\leq 4$,
which corresponds to  probabilities down to a level of $2\times 10^{-5}$  for the limiting
Maxwellian pdf.
Similarly,  if $g_{\text{tr}}=8$ then $\mathcal E_{\operatorname{rel}} < 10^{-1}$ 
for $v\leq 6$, corresponding to  probabilities down to $10^{-9}$. 
In this manner, the results in Fig.~\ref{fig:RelErrorContours} can be used to select a
value of $g_{\text{tr}}$ that is small but that nevertheless
guarantees a desired accuracy for the approximation  $Q^{\text{tr}}\approx Q$.
The advantage of choosing smaller values for $g_{\text{tr}}$
is that we can then  choose smaller values for $L$ and $N$, thereby
reducing  the computational cost, which is $\mathcal O(N^6)$.

In Fig.~\ref{fig:RelErrorContours} (right) we plot slices of 
$\mathcal E_{\operatorname{rel}}(v,g_{\text{tr}})$ for three values of $v$.
The colored solid curves with symbols show the results obtained with
Theorem~\ref{thmerror}, while the black dashed curves show the corresponding results obtained
using \eqref{corerror}. Even with $k=0.5$,
the asymptotic formulae \eqref{eq:Asymp3} agree extremely well 
with \eqref{Error_Estimate}, except when  $g_{\text{tr}}$ is slightly larger than $v$. 
These plots confirm that once ${g_{\text{tr}}} > v$,  
$Q^{{\text{tr}}}(\mathbf v)\to Q(\mathbf v)$ exponentially fast as ${g_{\text{tr}}}\to\infty$.
However, the pointwise nature of the convergence is obvious in the plots. 

\begin{figure}[t!]
\centering
\includegraphics[width=0.49\textwidth]{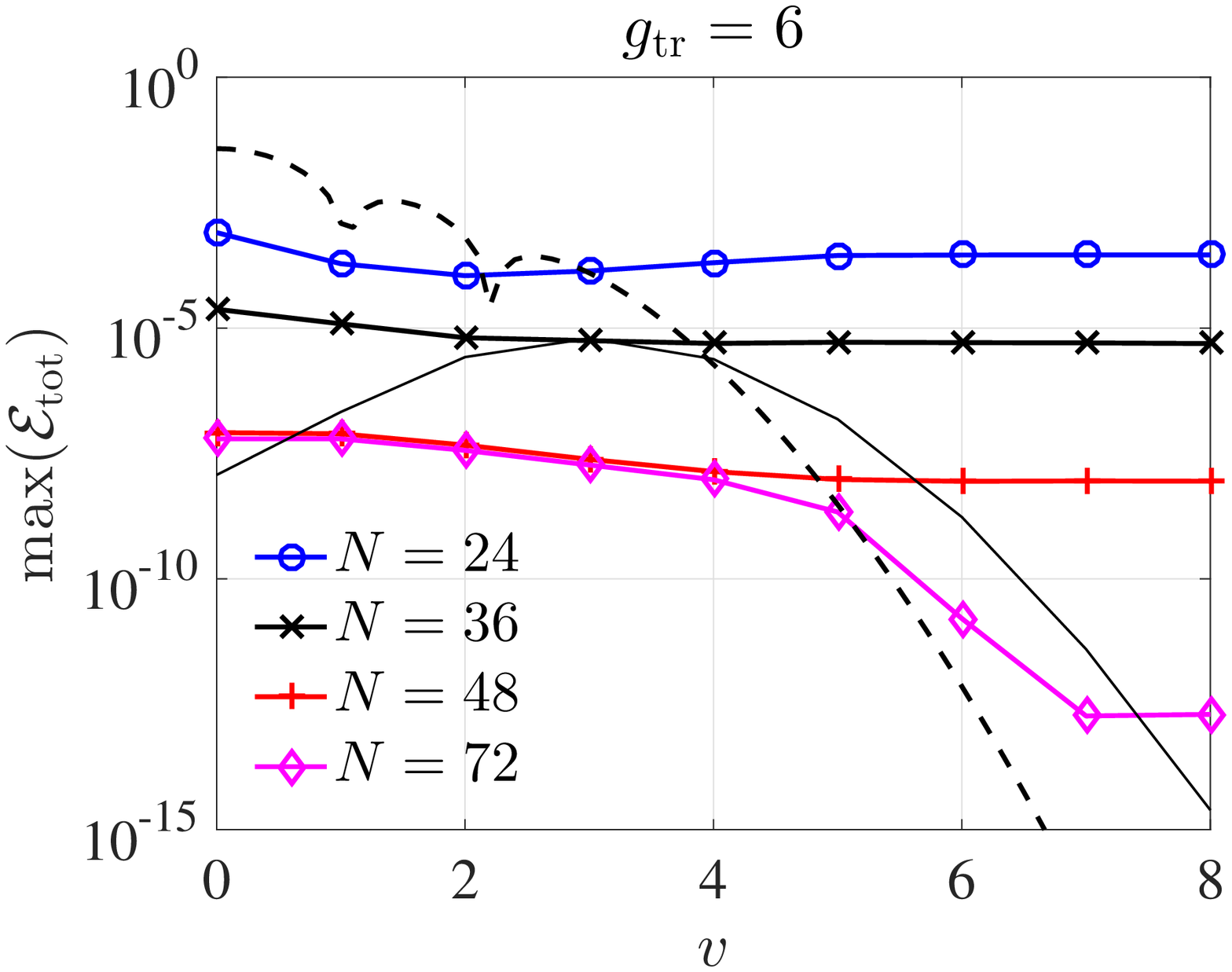} 
\includegraphics[width=0.49\textwidth]{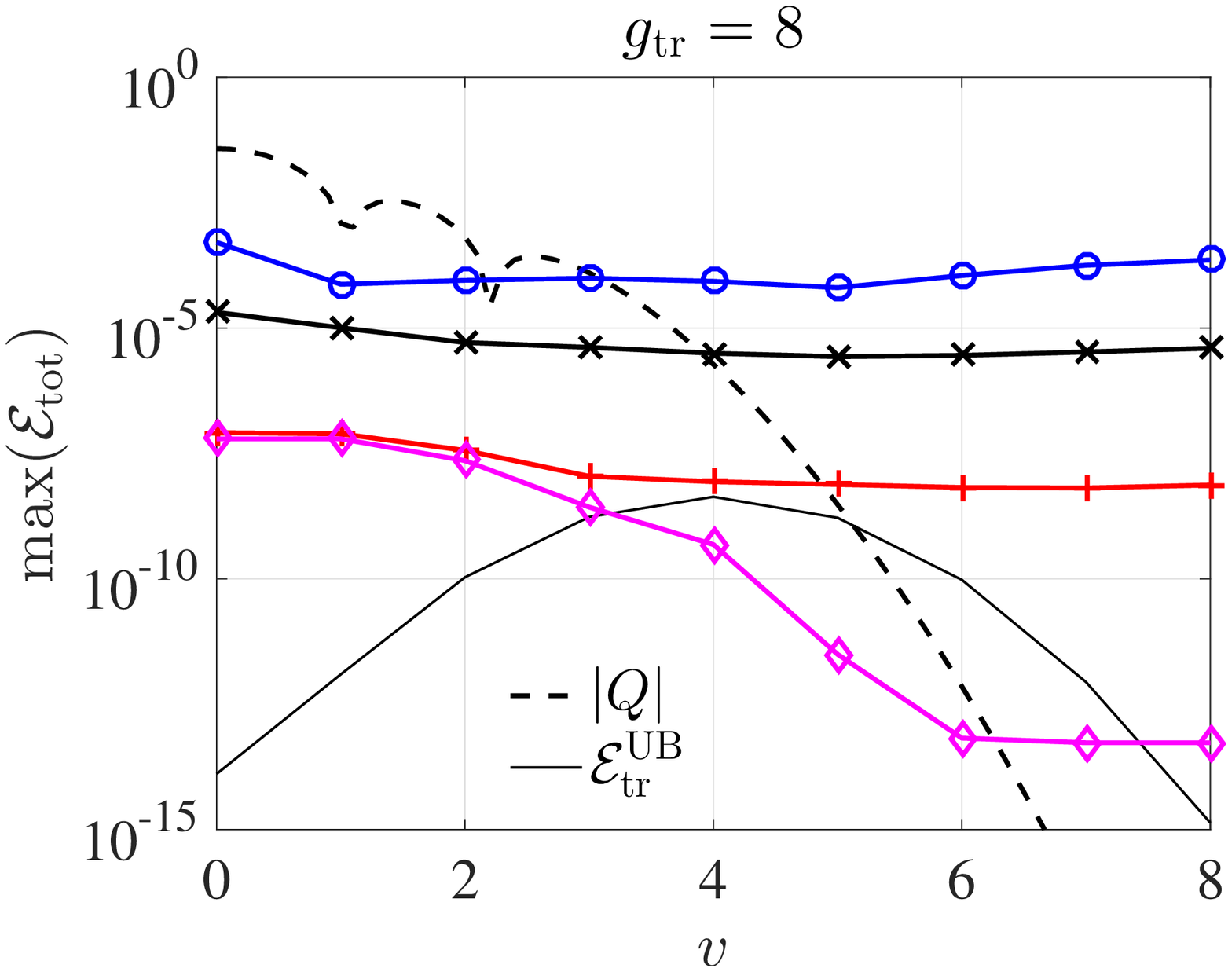} \\
\includegraphics[width=0.49\textwidth]{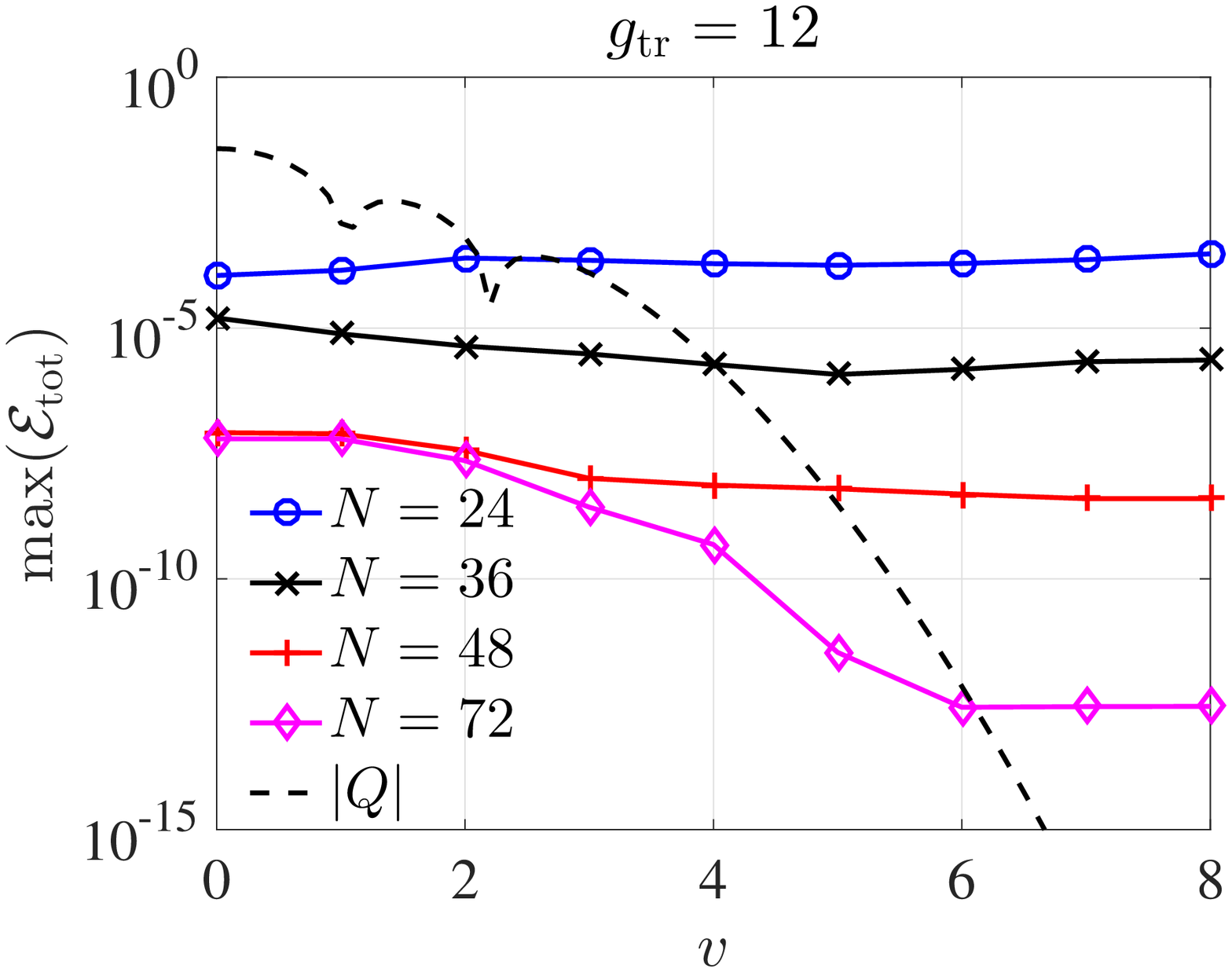} 
\includegraphics[width=0.49\textwidth]{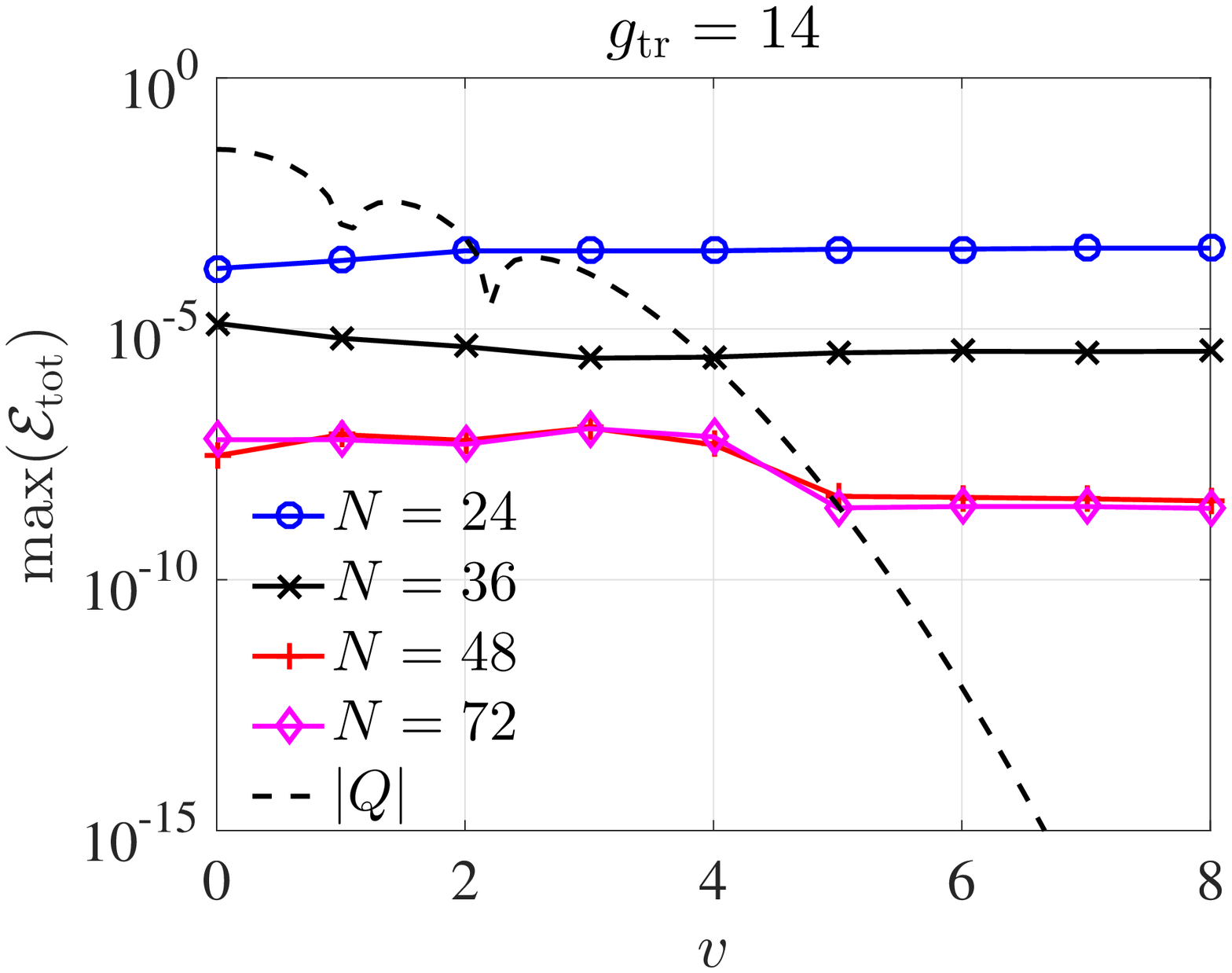} 
\caption{Maximum of the total error,
$\mathcal E_{\operatorname{tot}}$ in \eqref{TotErr}, as a function of speed $v$, 
for several different choices of $g_{\operatorname{tr}}$ and $N$.
These results are for the BKW solution.
Here, the maximum is taken over all $\mathbf v$ with $|\mathbf v|=v$.
The solid colored curves with symbols show plots of 
$\operatorname{max}(\mathcal E_{\operatorname{tot}})$
for the values of $N$ shown in the legend. 
The black dashed curve shows the collision operator, $Q$, obtained from \eqref{eq:HBE} and \eqref{eq:BKW}, and the solid black curve shows the upper bound, 
$\mathcal E_{\operatorname{tr}}^{\operatorname{UB}}$ in \eqref{anarelerror}, 
for the truncation error. 
}
 \label{fig:BKWAbsErrUBQ}
\end{figure}

In Fig.~\ref{fig:BKWAbsErrUBQ} we assess the accuracy of the numerical
computation of $Q^{\operatorname{tr}}$ by plotting the 
maximum of the total error,
$\mathcal E_{\operatorname{tot}}$ in \eqref{TotErr}, as a function of speed $v$, 
for several different choice of 
$g_{\operatorname{tr}}$ and $N$.
Here, the maximum is taken over all $\mathbf v$ with $|\mathbf v|=v$.
In the top left panel, we show the results with
$g_{\operatorname{tr}}=6$. Using solid colored curves with symbols 
we plot $\operatorname{max}(\mathcal E_{\operatorname{tot}})$
for $N=24$ (blue curve with circles), $N=36$ (black curve with crosses),
$N=48$ (red curve with pluses), and $N=72$ (magenta curve with diamonds).
We also plot the collision operator, $Q$, obtained analytically from 
\eqref{eq:HBE} and \eqref{eq:BKW}
(dashed black curve) and the upper bound, 
$\mathcal E_{\operatorname{tr}}^{\operatorname{UB}}$ in \eqref{anarelerror}, 
for the truncation error  (solid black curve). 
Because these two curves intersect
at $v=4$, we can only be guaranteed that 
$|Q^{\operatorname{tr}} - Q| < |Q|$ for $v\leq 4$.
For each $N$, the numerically computed
collision operator, $Q^{\operatorname{NC}}$, is an accurate approximation to 
$Q$ in the interval where the solid curve with symbols lies below the 
black dashed curve. The reason this interval extends past $v=4$ for
$N\geq 48$ is that the upper bound $\mathcal E_{\operatorname{tr}}^{\operatorname{UB}}$  for $|Q^{\operatorname{tr}} - Q| $ is not optimal. 
In the top right panel, we show the corresponding results with
$g_{\operatorname{tr}}=8$. 
Because the solid black and dashed black curves intersect at $v=5$, we are
guaranteed that if we choose $N$ to be sufficiently large, then the
solution will be accurate out to at least $v=5$. Clearly, the choice $N=48$
is not large enough. However, if we choose $N=72$, corresponding to a 12-fold
increase in the computational time, then the solution is accurate out to $v=6$. 
%%
\begin{comment}
Because the intersection point of the red curve with pluses 
and the black dashed  curve occurs
at a slightly smaller value of $v$ than the
intersection point of the two black  curves,
we know that if we choose $N$ to be sufficiently larger than $48$
the solution will be accurate out to at least $v=5$.
Indeed,  with $N=72$ and an increase in computational time by a factor of 12, 
we obtain agreement out to $v=6$.
\end{comment}
%%
In the bottom left panel for which $g_{\operatorname{tr}}=12$,
the solid black curve is not visible since 
$\mathcal E_{\operatorname{tr}}^{\operatorname{UB}} < 10^{-15}$.
However, there is no change in the
$\mathcal E_{\operatorname{tot}}$-curves compared to the case that
$g_{\operatorname{tr}}=8$, since the total error is dominated by the
error in the numerical computation of $Q^{\operatorname{tr}}$.
Finally, in the bottom right panel with $g_{\operatorname{tr}}=14$,
we see that there is no advantage to increasing $N$ from 48 to 72
since that does not decrease $\Delta \xi$ and the convolution
weighting function, ${\widehat G}^{\operatorname{tr}}$, now
oscillates too rapidly. This last result is in accord with the large
jump in the errors from $g_{\operatorname{tr}}=12$ to $g_{\operatorname{tr}}=16$
that we observed for the Maxwellian pdf in Table~\ref{table:MaxMaxwellianQBeforeLP}.

\begin{figure}[t!]
\centering
\includegraphics[width=0.49\textwidth]{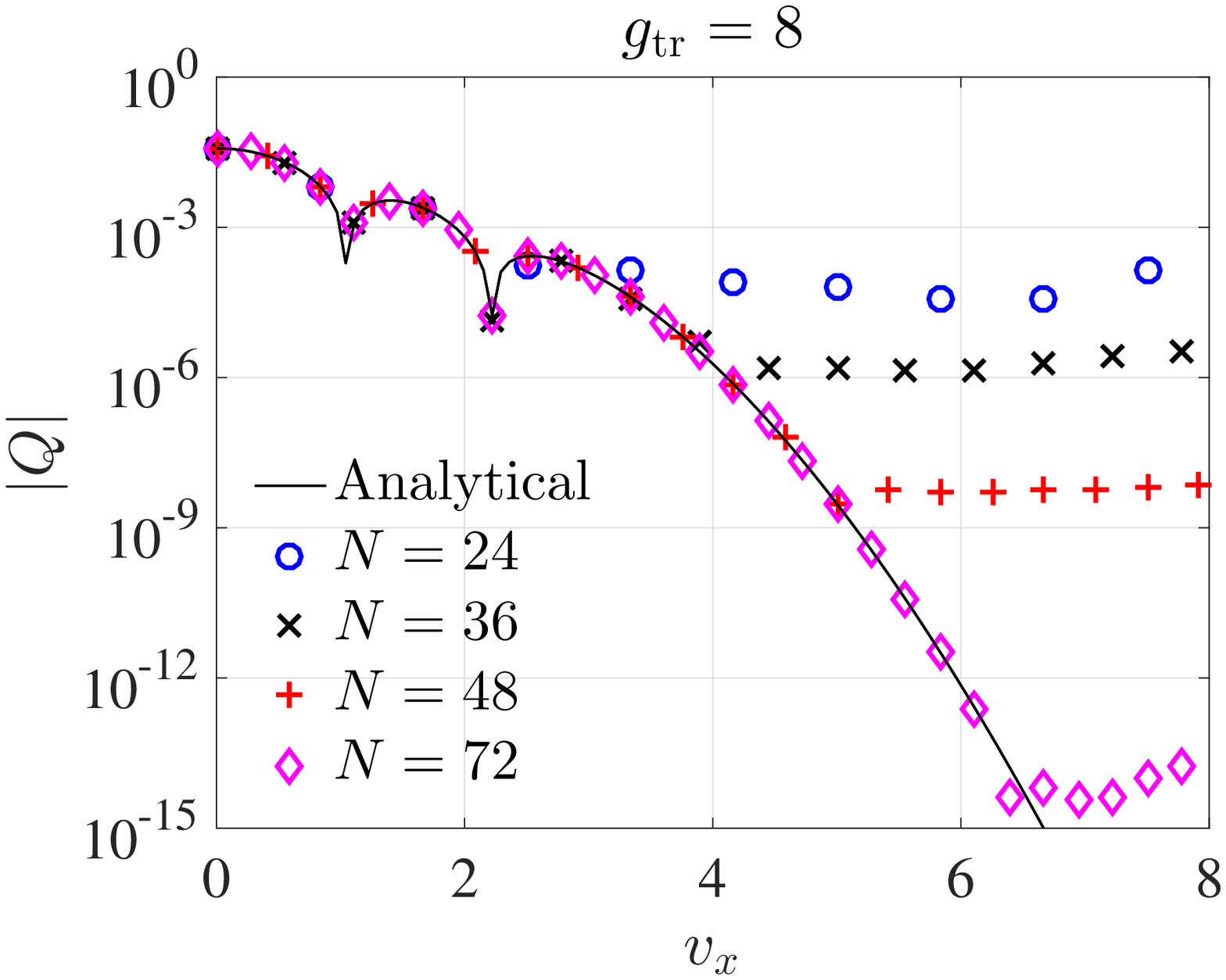}
\includegraphics[width=0.49\textwidth]{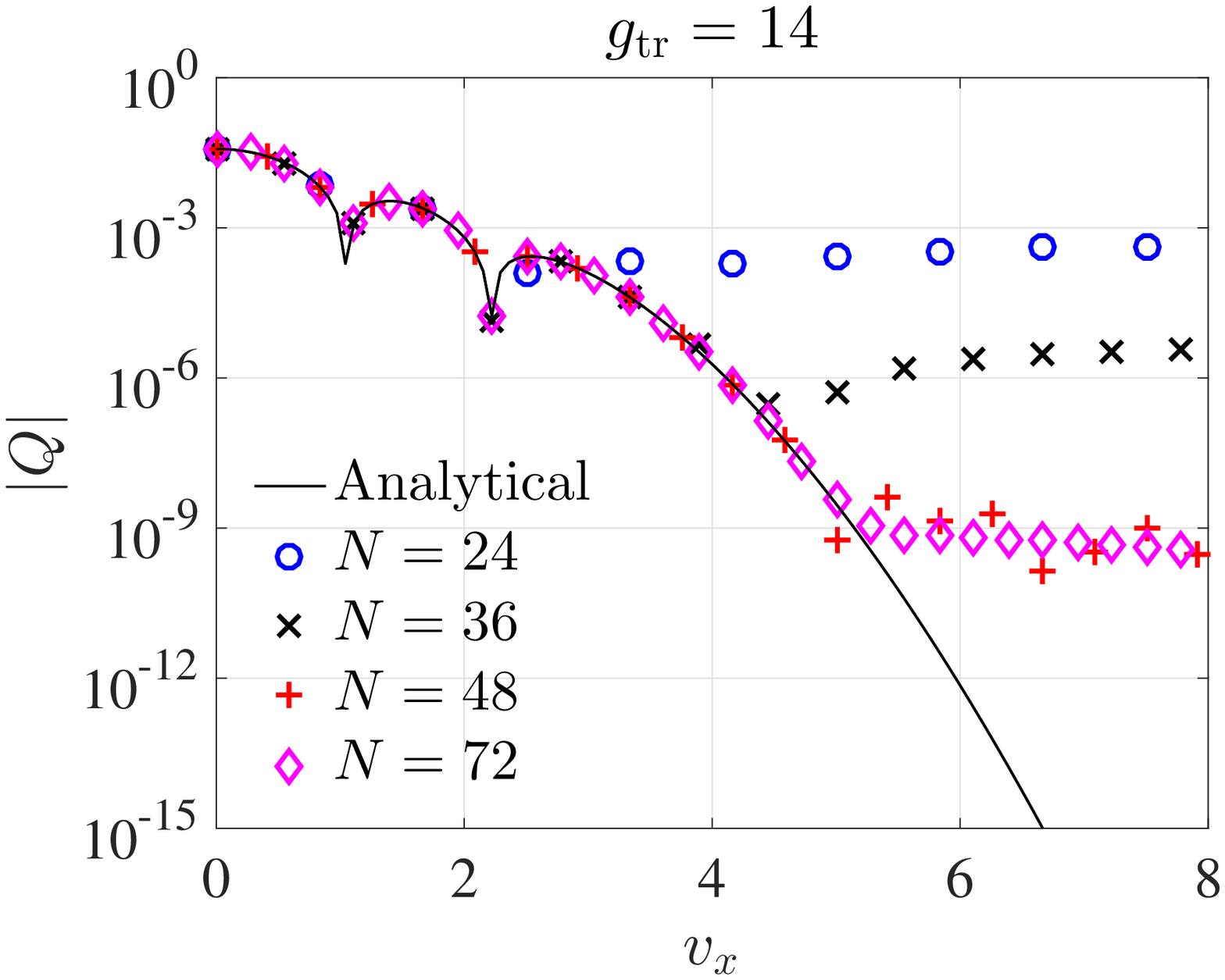}\\
\includegraphics[width=0.49\textwidth]{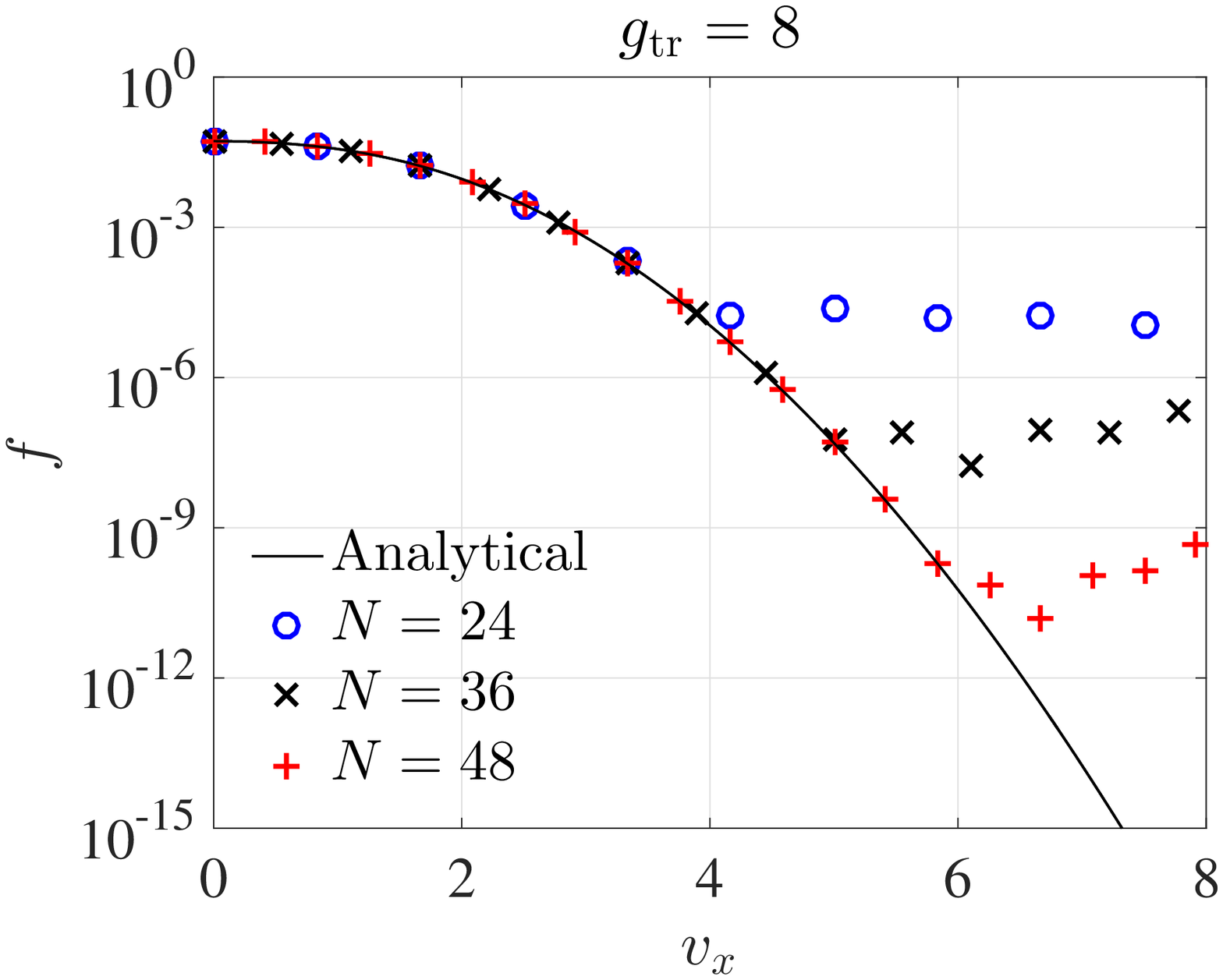}
\includegraphics[width=0.49\textwidth]{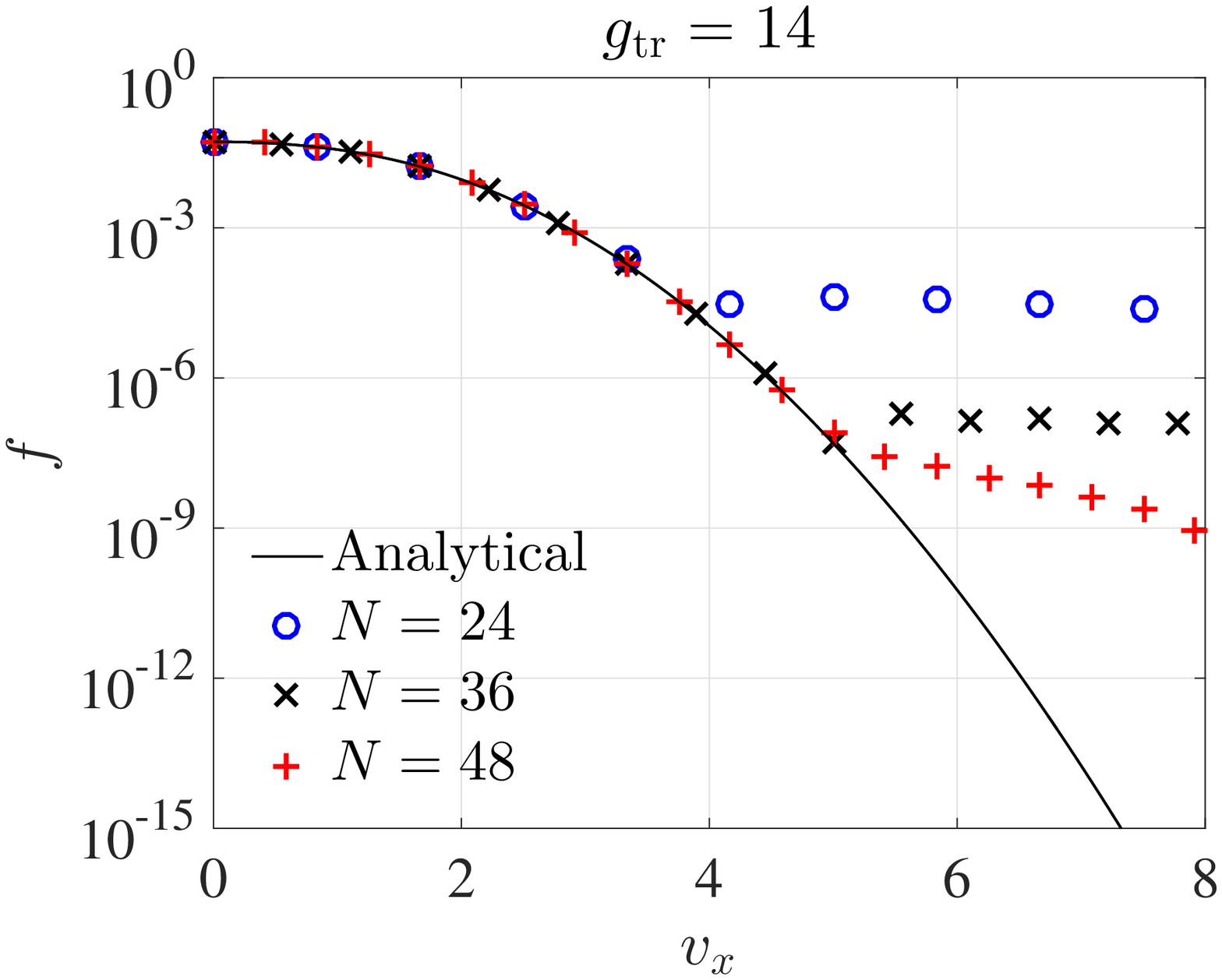}
\caption{Collision operator at $t=5.5$ (top row) and velocity pdf at $t=9$
(bottom row) for the BKW solution. The numerical results were obtained
with $g_\text{tr}=8$ (left column) and $g_\text{tr}=14$ (right column) for the
values of $N$ shown in the legends. The analytical solutions are
shown with the black solid curves.}
 \label{fig:BKWQ}
\end{figure}

In the top row of Fig.~\ref{fig:BKWQ}, we plot $|Q|$ as a function of $v_x$
at $(v_y,v_z)=(0,0)$. We show the numerical results obtained 
with $g_\text{tr}=8$ (left) and $g_\text{tr}=14$ (right) 
for the values of $N$ shown in the legend. 
We also show the analytical result obtained  from \eqref{eq:HBE}
and  \eqref{eq:BKW}
with a black solid curve. The cusps correspond to the values of $v_x$ for which
$Q=0$. 
On a linear scale (not shown), we obtain excellent agreement for all values of $N$.
With $N=48$, we obtain excellent agreement down to the level of less than
$10^{-8}$, and with $N=72$ down to $10^{-14}$. 
The results with $g_\text{tr}=6$, 10, and 12 (not shown) are only slightly worse
than with $g_\text{tr}=8$.  However, just as in Fig.~\ref{fig:BKWAbsErrUBQ},
with $g_\text{tr}=14$ we cannot
reduce the error level  below $10^{-9}$. 
In the bottom row, we plot the velocity pdf at $t=9$, using the same format
as in the top row. 
 For these results
we solved \eqref{eq:HBE} using Euler's method with a time step of $\Delta t = 0.05$. 
However, we did not perform the computation with
$N=72$ as the computational cost was prohibitive.

\subsection{A cylindrically symmetric initial condition}\label{Sec:CS}

In this subsection, we apply the  spectral method to solve the 
homogeneous Boltzmann equation~\eqref{eq:HBE} in the case that the  initial condition
is the  cylindrically symmetric velocity pdf,
\begin{equation}
{f}_{\text{cyl}}(0,\boldsymbol{v}) \,\,=\,\, A 
 \left[ \tfrac{5K-3}{K} +\tfrac{1-K}{K^2}\,(c^2 v_x^2+c^2v_y^2+v_z^2)\right]
\exp\left[-(c^2 v_x^2+c^2v_y^2+v_z^2)/2K\right],
\label{eq:CS}
\end{equation}
obtained by dilating the BKW initial condition~\eqref{eq:BKW} by a factor, $1/c$, 
in the $v_x$ and $v_y$-dimensions. 
As in Section~\ref{Sec:BKW}, we choose $K=1-e^{-5.5/6}$. 
We choose the dilation constant to be $c=2$ and we choose the
constant, $A$, so that the pdf integrates to 1.

\begin{figure}[t!]
\centering
\includegraphics[width=0.32\textwidth]{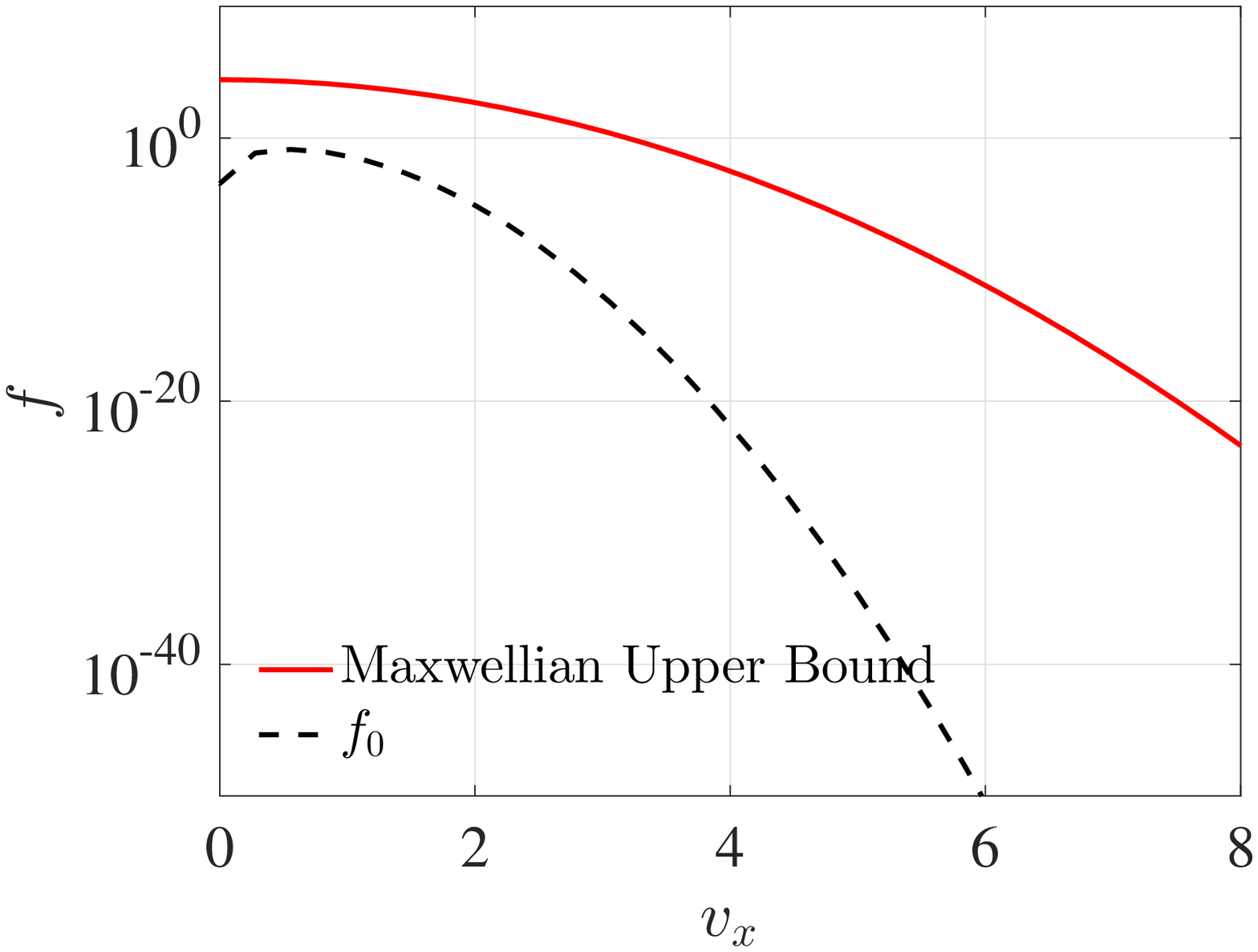}
\includegraphics[width=0.32\textwidth]{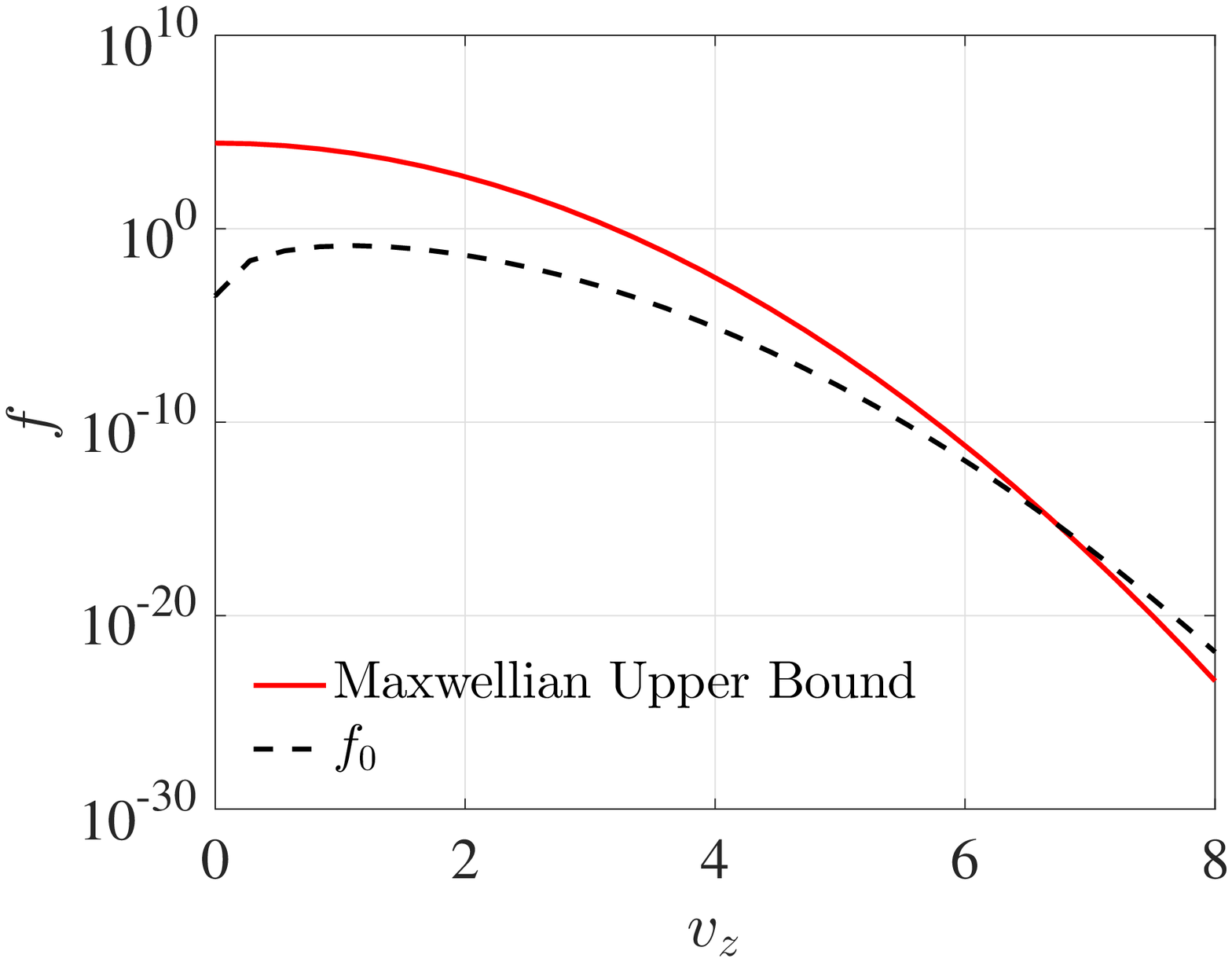}
\includegraphics[width=0.32\textwidth]{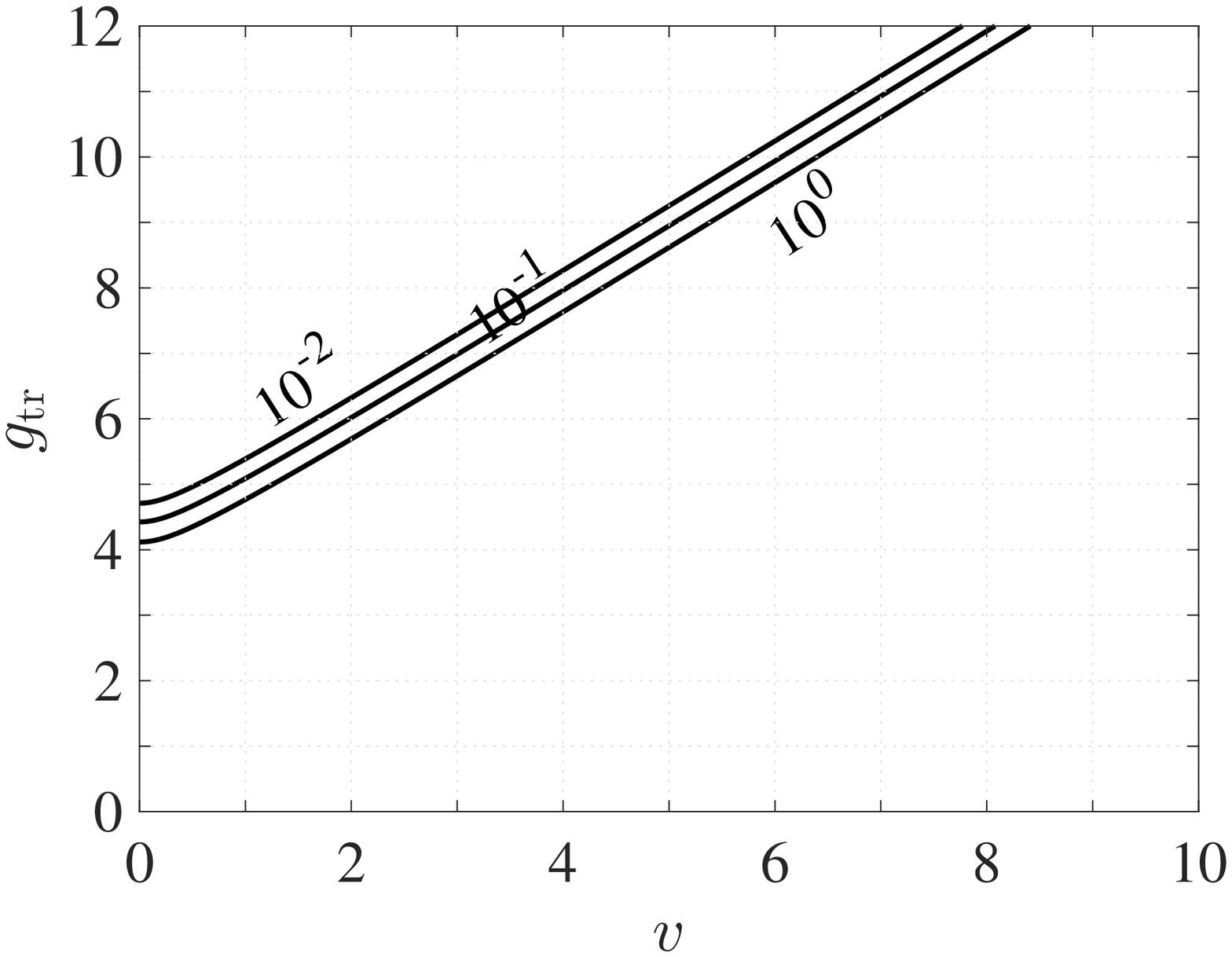}
\caption{Left and Middle: Log-scale plot of the cylindrically symmetric
initial velocity pdf \eqref{eq:CS}
(dashed black curve) and the Maxwellian upper bound  (solid red curve)
as functions of $v_x$ when $(v_y,v_z)=(0,0)$ (left) and $v_z$ 
when $(v_x,v_y)=(0,0)$ (middle).
Right: Contour plot of the upper bound, $\mathcal E_{\operatorname{rel}}$, for the relative error in the truncation of the collision operator given by \eqref{Error_Estimate},
as a function of speed, $v$, and truncation parameter, $g_{\rm tr}$.
}
 \label{fig:CylSymmContour}
\end{figure}

\begin{figure}[h!]
\centering
\includegraphics[width=0.49\textwidth]{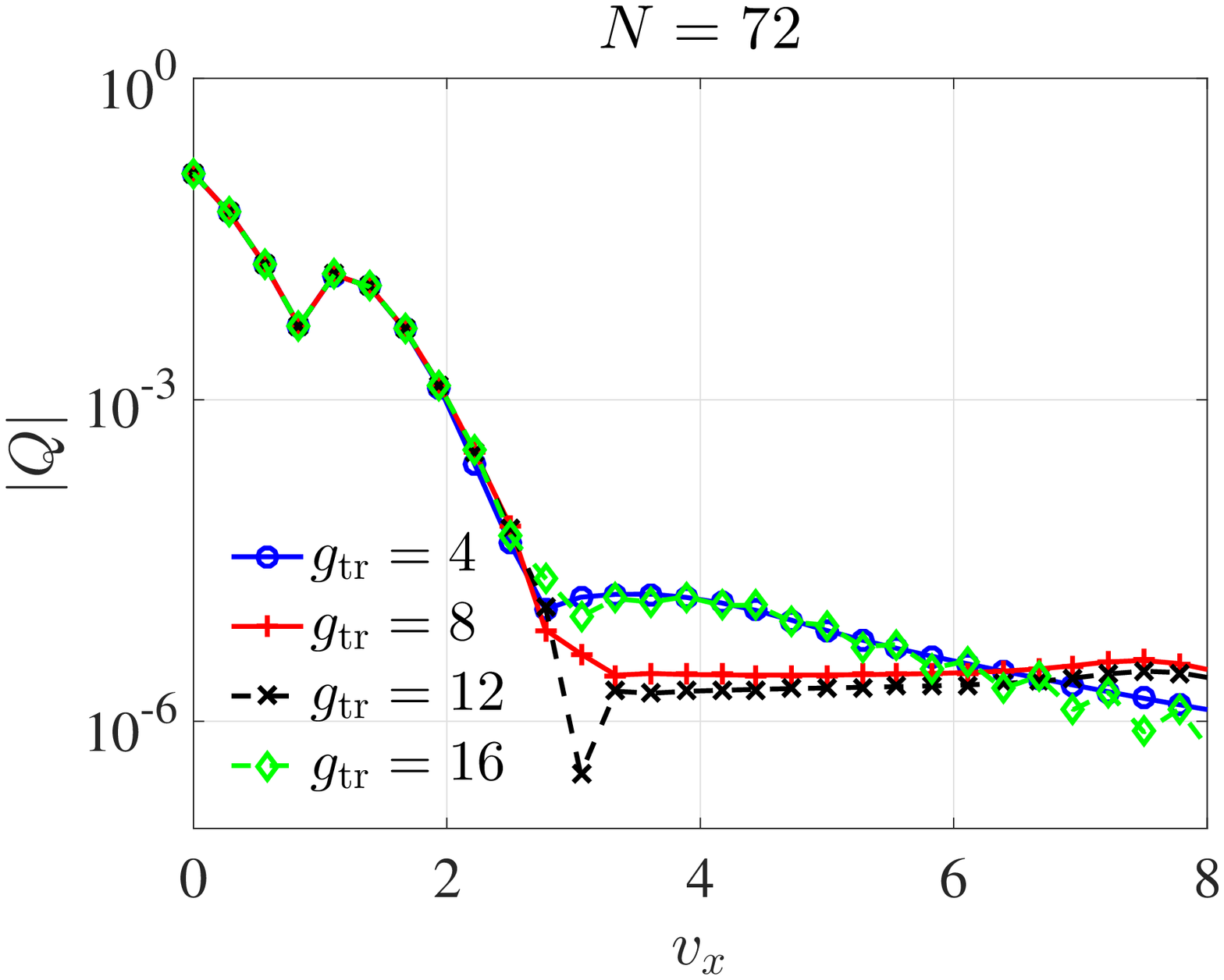}
\includegraphics[width=0.49\textwidth]{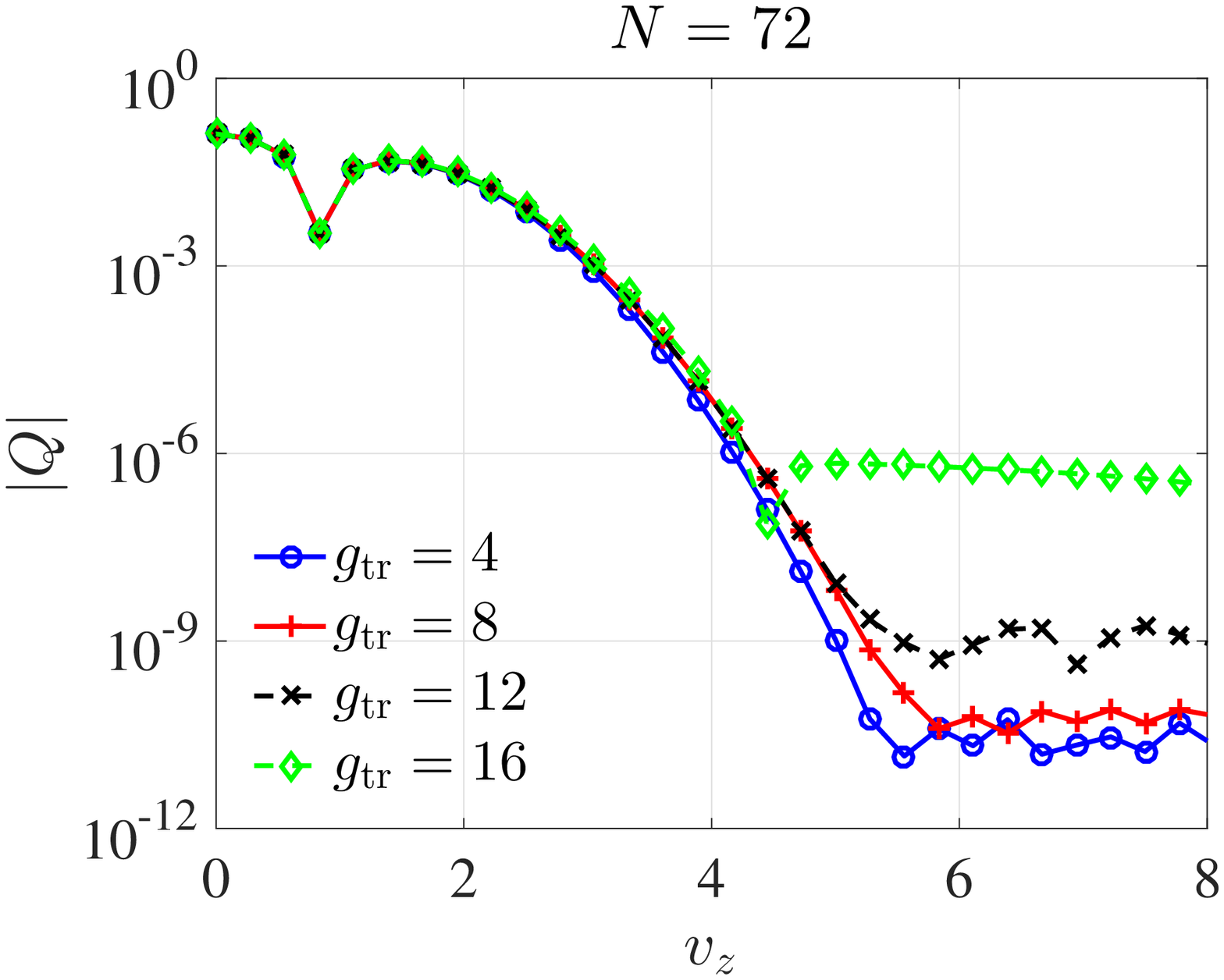}\\
\includegraphics[width=0.49\textwidth]{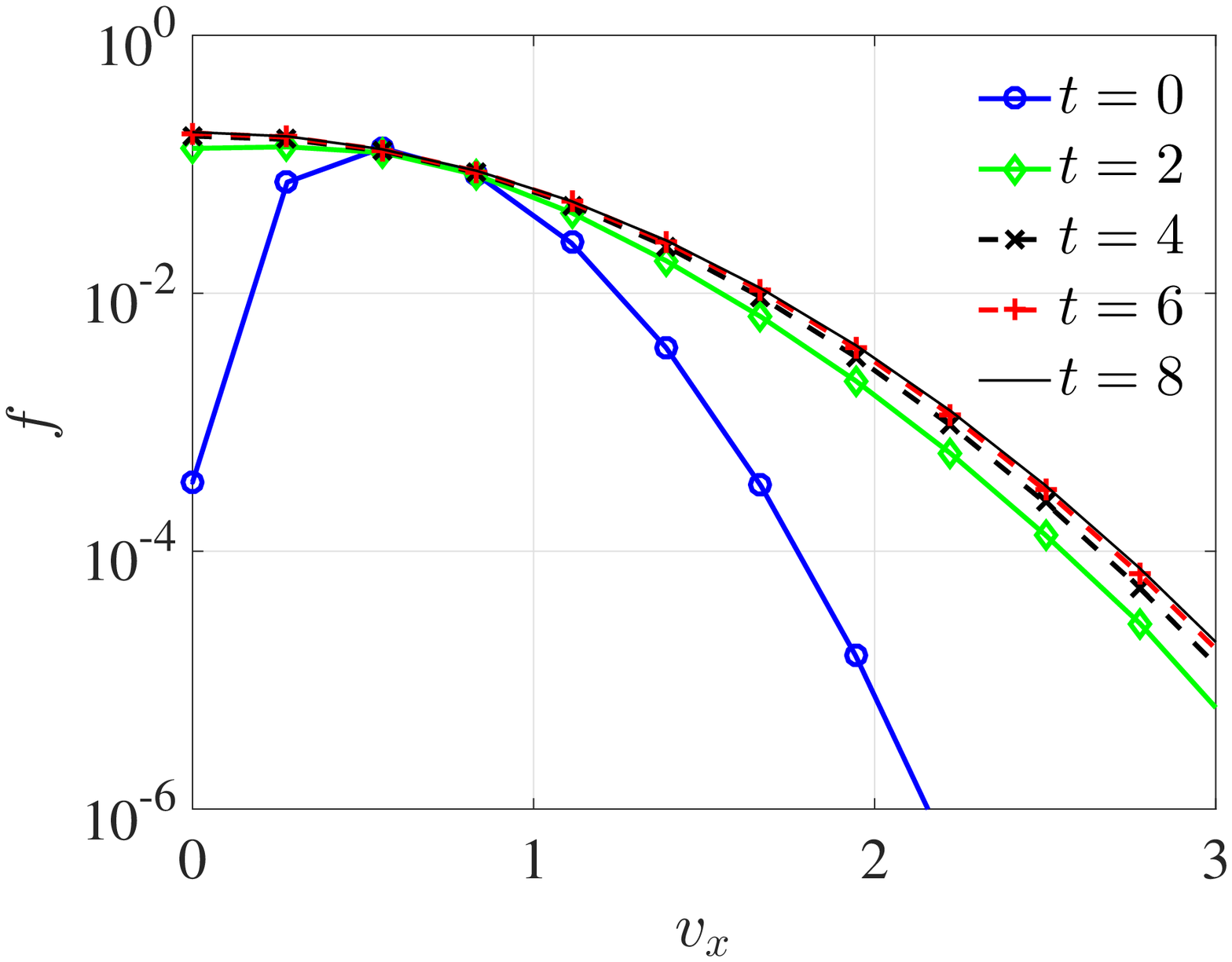}
\includegraphics[width=0.49\textwidth]{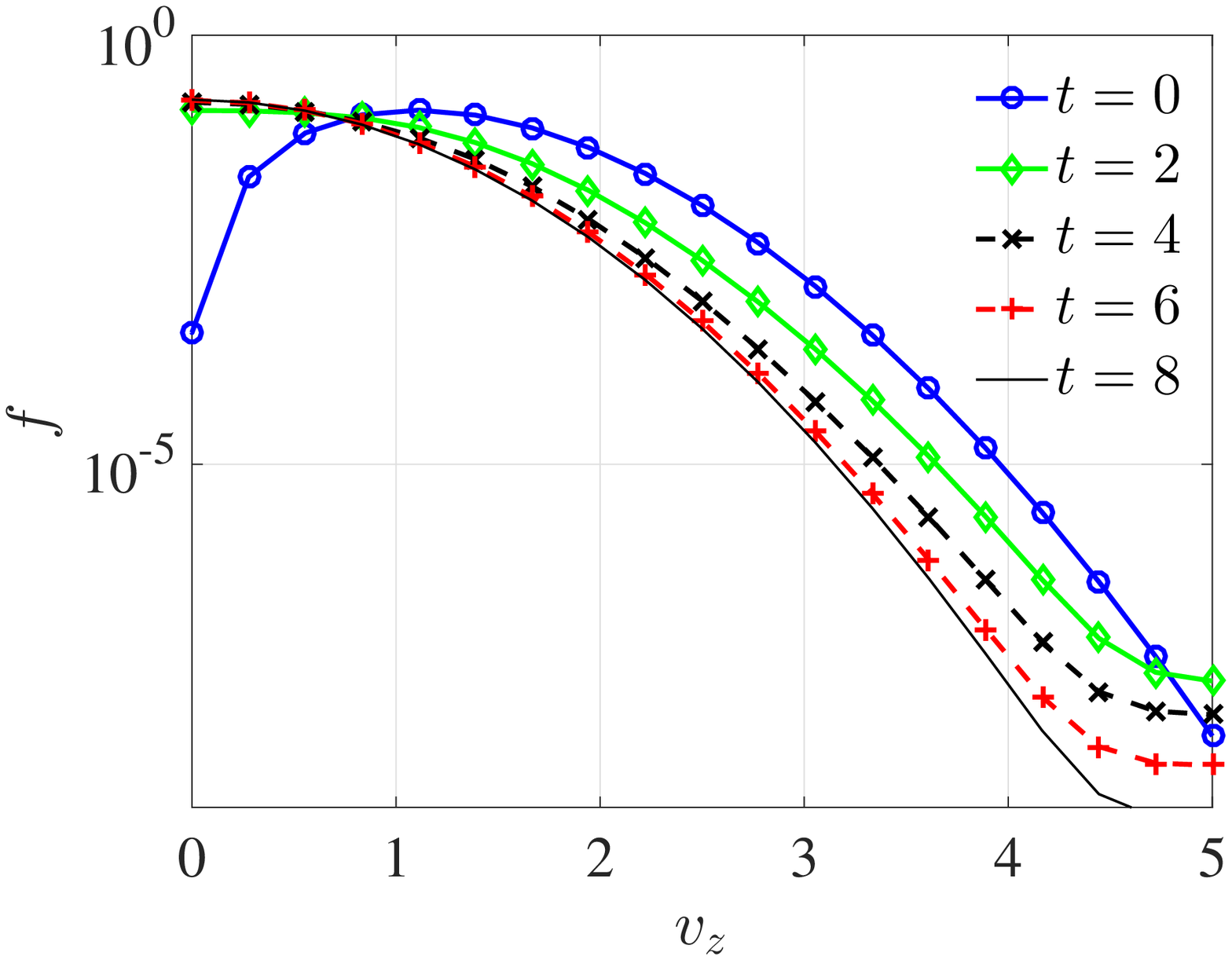}\\
\caption{Top row: Collision operator at $t=5.5$ 
as a function of $v_x$ (left) and $v_z$ (right) 
for the cylindrically symmetric initial velocity pdf~\eqref{eq:CS}
with $c=2$ and $N=72$, for the values of $g_\text{tr}$ shown in the legends.
Bottom row: The corresponding velocity pdfs at the times shown in the legends,
computed using $g_\text{tr}=12$. 
}
 \label{fig:CylSymmQVaryGtrAndTimeSeries}
\end{figure}

In Fig.~\ref{fig:CylSymmContour}, we plot the initial pdf, $f_{\operatorname{cyl}}$,
and the Maxwellian upper  bound obtained using Method~I, which  gives $k=0.9$ and 
$c=1.5\times 10^5$. 
We plot these pdfs as a function of
both $v_x$ when $(v_y,v_z)=(0,0)$ (left) and $v_z$ when $(v_x,v_y)=(0,0)$ (middle). 
In the right panel, 
we show a contour plot of the bound, $\mathcal E_{\operatorname{rel}}$, 
for the  relative error in the truncation of the collision operator 
given by~\eqref{Error_Estimate}, as a function of  $v$ and $g_{\text{tr}}$. 
We observe that the contours are translated up by about 2 compared
to the ones in Fig.~\ref{fig:RelErrorContours}.
The contour plot shows that if we choose $g_{\text{tr}}=10$ then
$\mathcal E_{\operatorname{rel}} < 10^{-1}$ for $v\leq 6$.

In the top row of Fig.~\ref{fig:CylSymmQVaryGtrAndTimeSeries}, we
plot the numerical collision operator as a function of $v_x$ (left) and 
$v_z$ (right). These results were obtained using $N=72$ for
the values of $g_{\text{tr}}$ shown in the legend. 
The results for $g_{\text{tr}}=4$ and 16 are less accurate than
those for  $g_{\text{tr}}=8$ and 12. 
The results for $g_{\text{tr}}=4$ is consistent
with the contour plot in Fig.~\ref{fig:CylSymmContour}, which shows that
the truncation error is too large when $g_{\text{tr}}=4$. 
Because of the more rapid decay of the cylindrically symmetric initial
condition in the $v_x$-direction, the error in the numerical computation
of $Q^{\operatorname{tr}}$ is the dominant source of error
with $g_{\text{tr}}=8$ and 12.
In the bottom row of Fig.~\ref{fig:CylSymmQVaryGtrAndTimeSeries}, we
plot the evolution of the velocity pdf over the time interval
$[0,8]$. 
These results were obtained with $g_{\text{tr}}=12$,
using the Adams-Bashforth method~\eqref{Adams}
with $\Delta t = 0.125$. 
We verified that the
number density, momentum, and energy are preserved up to round-off
error. 
The pdf at $t=8$ agrees well with the equilibrium Maxwellian pdf (not shown)
over the range of probability values in the plots. 
The plots show the rates at which the velocity pdf converges to the
equilibrium pdf  in the different velocity dimensions.

\subsection{Mixture of Maxwellians initial conditions}\label{Sec:MixMax}
For the next two examples, we suppose that the initial velocity pdf
is a mixture of Maxwellian pdfs  of the form
\begin{equation}\label{eq:MixMax}
f_{\operatorname{mix}}(\boldsymbol{v}) \,\,=\,\, 
\omega f(\boldsymbol{v}-\boldsymbol{v}_1,T_1)
\,\,+ \,\,(1-\omega) f(\boldsymbol{v}-\boldsymbol{v}_2,T_2),
\end{equation}
where $f(\boldsymbol{v},T) = (2\pi T)^{-3/2} \exp(-v^2/2T)$.

\begin{figure}[t!]
\centering
\includegraphics[width=0.32\textwidth]{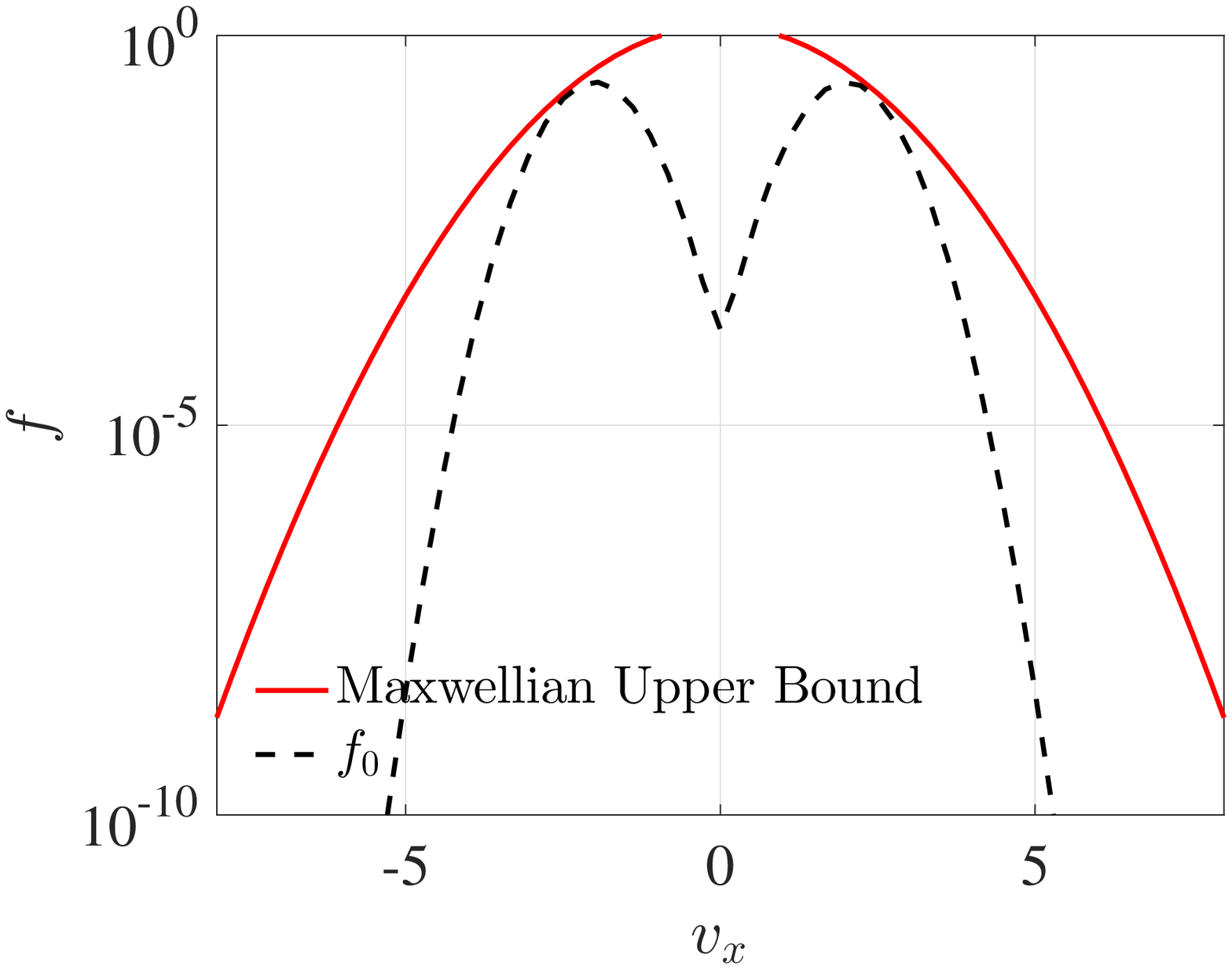}
\includegraphics[width=0.32\textwidth]{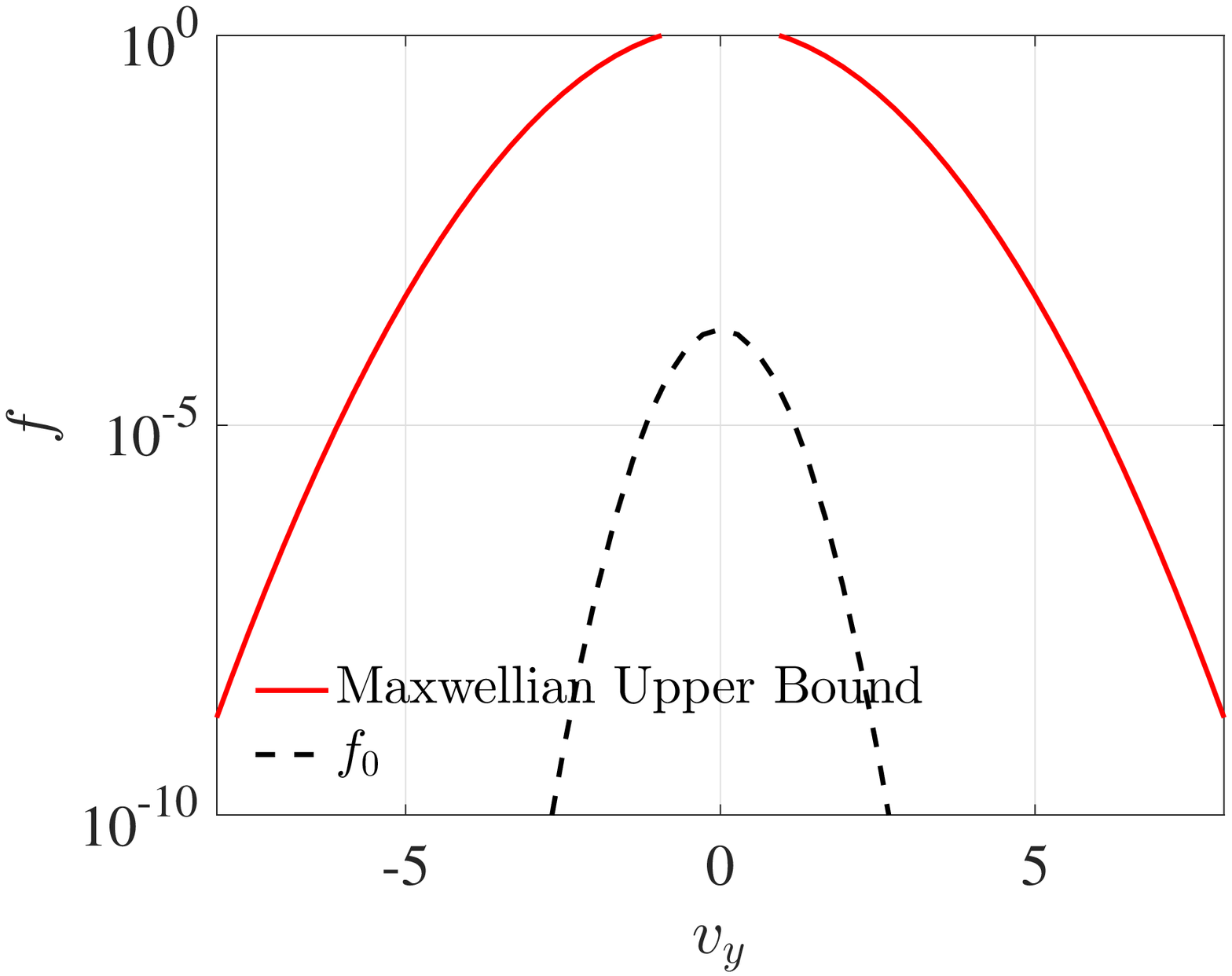}
\includegraphics[width=0.32\textwidth]{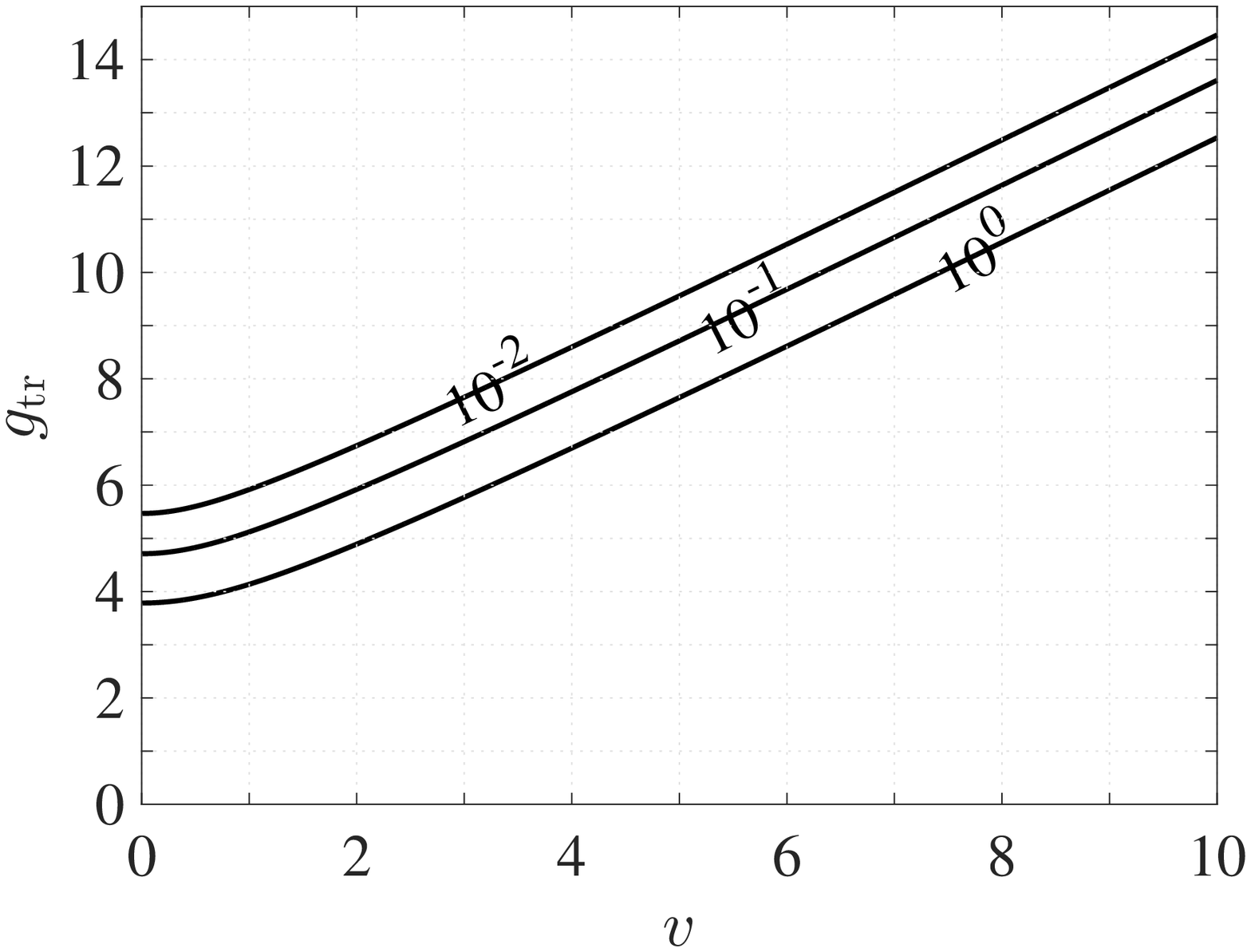}
\caption{
Left and Middle: Log-scale plot of the mixture of Maxwellians 
initial velocity pdf~\eqref{eq:MixMax}
(dashed black curve) and the Maxwellian upper bound   (solid red curve) 
as functions of $v_x$ when $(v_y,v_z)=(0,0)$ (left) and $v_y$ 
when $(v_x,v_z)=(0,0)$ (right). 
The parameters in~\eqref{eq:MixMax} were chosen to be
$\omega=0.5$, $\boldsymbol{v}_1 = (2,0,0)$,
$\boldsymbol{v}_2=-\boldsymbol{v}_1$, and $T_1=T_2=0.25$. 
Right: Contour plot of the upper bound, $\mathcal E_{\operatorname{rel}}$, for the relative error in the truncation of the collision operator given by \eqref{Error_Estimate}.
%%Mix X2. Suggests $g_\text{tr}=10$. 
}
 \label{fig:MixMaxX2Contour}
\end{figure}

For our first example, we chose $\omega=0.5$, $\boldsymbol{v}_1 = (2,0,0)$,
$\boldsymbol{v}_2=-\boldsymbol{v}_1$, and $T_1=T_2=0.25$. 
In Fig.~\ref{fig:MixMaxX2Contour}, we plot the initial pdf, $f_{\operatorname{mix}}$,
and the Maxwellian upper  bound obtained using Method~I, which  gives $k=0.32$ and 
$c=1.1$.  
We plot these pdfs as a function of
both $v_x$ when $(v_y,v_z)=(0,0)$ (left) and $v_y$ when $(v_x,v_z)=(0,0)$ (middle). 
In the right panel, 
we show a contour plot of the bound, $\mathcal E_{\operatorname{rel}}$, 
for the  relative error in the truncation of the collision operator 
given by~\eqref{Error_Estimate}, as a function of  $v$ and $g_{\text{tr}}$. 
Guided by this contour plot, for the computation of the 
velocity pdf we chose $g_{\text{tr}}=10$
to ensure that $\mathcal E_{\operatorname{rel}} < 10^{-1}$ for $v\leq 6$.

\begin{figure}[h!]
\centering
\includegraphics[width=0.49\textwidth]{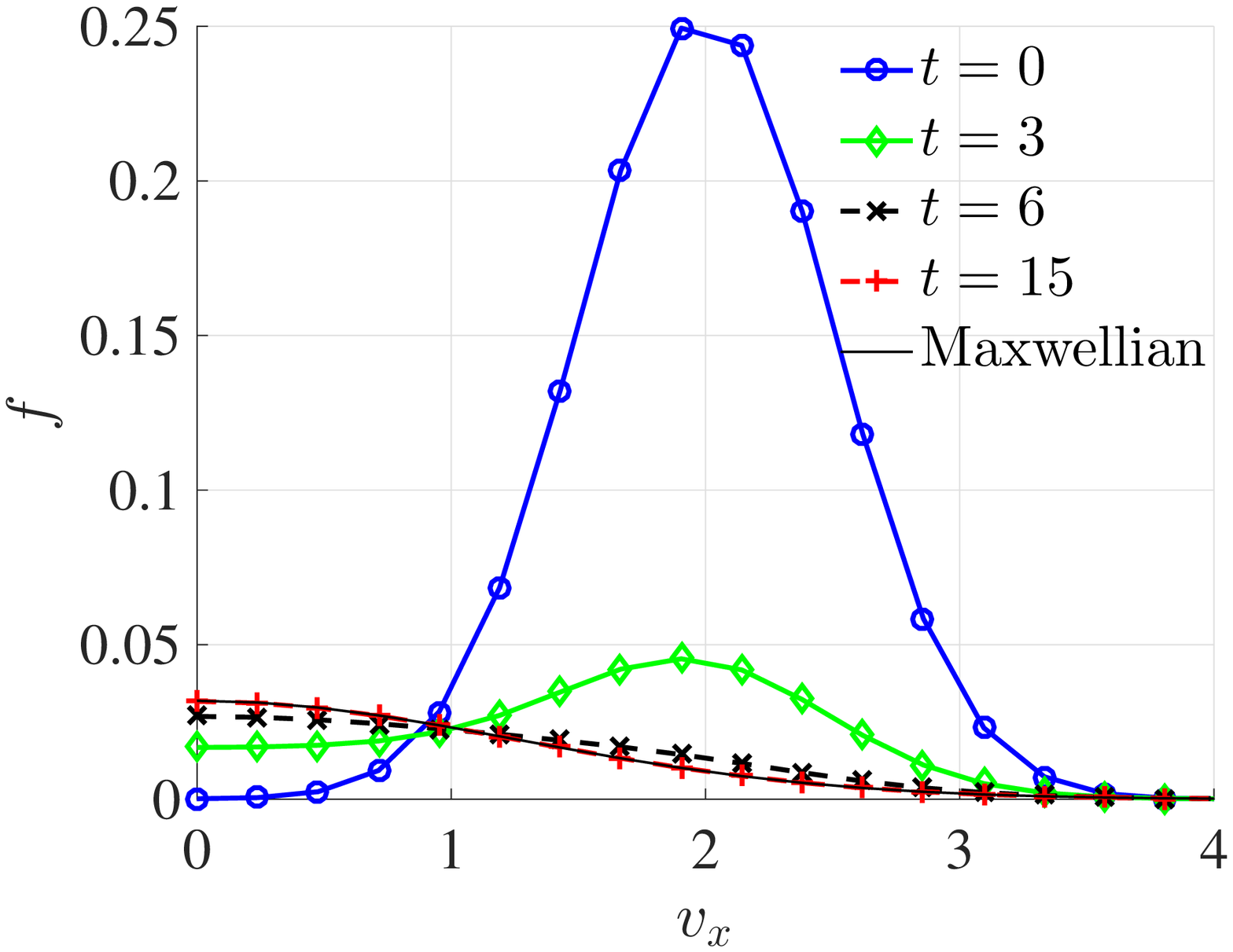}
\includegraphics[width=0.49\textwidth]{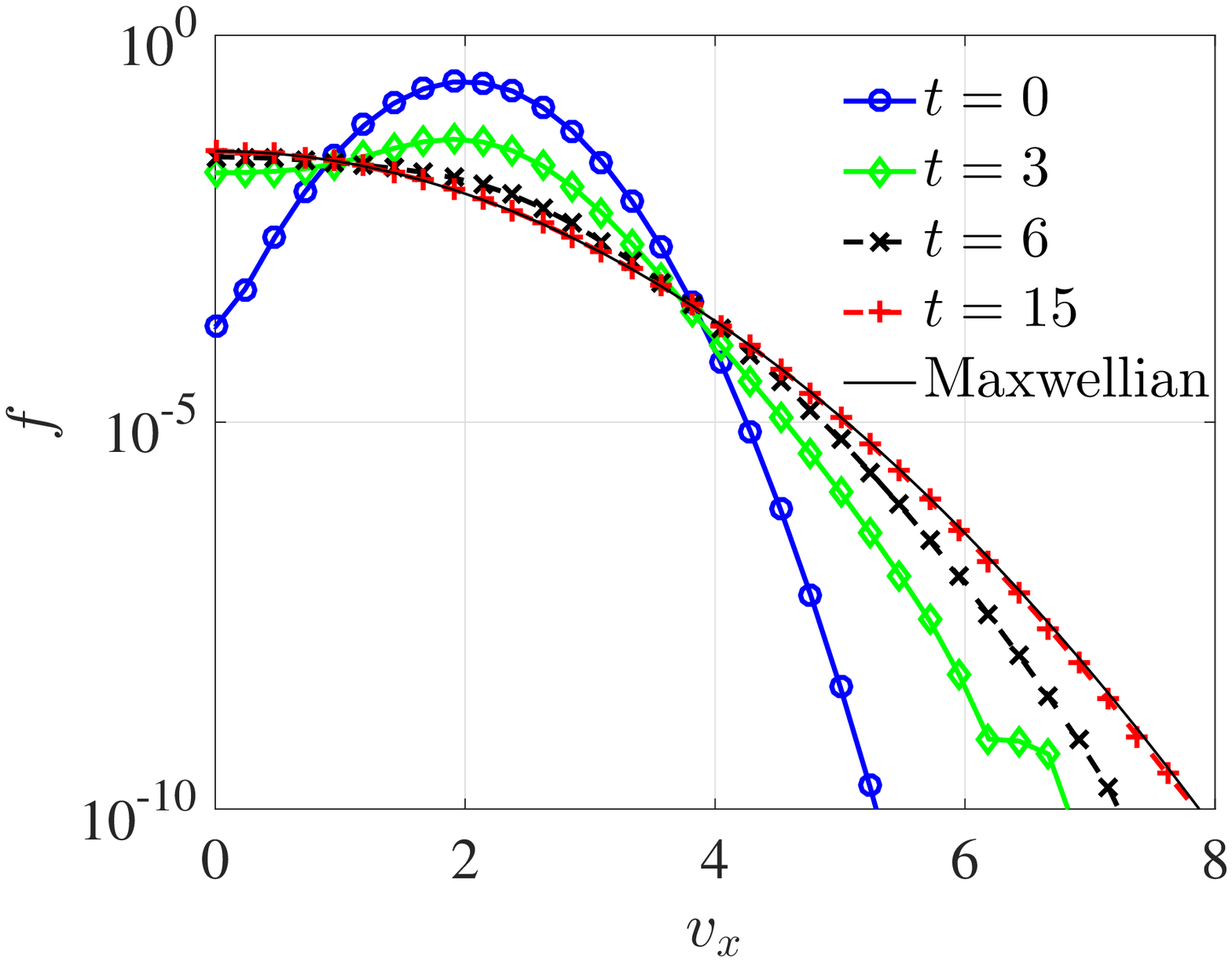}\\
\includegraphics[width=0.49\textwidth]{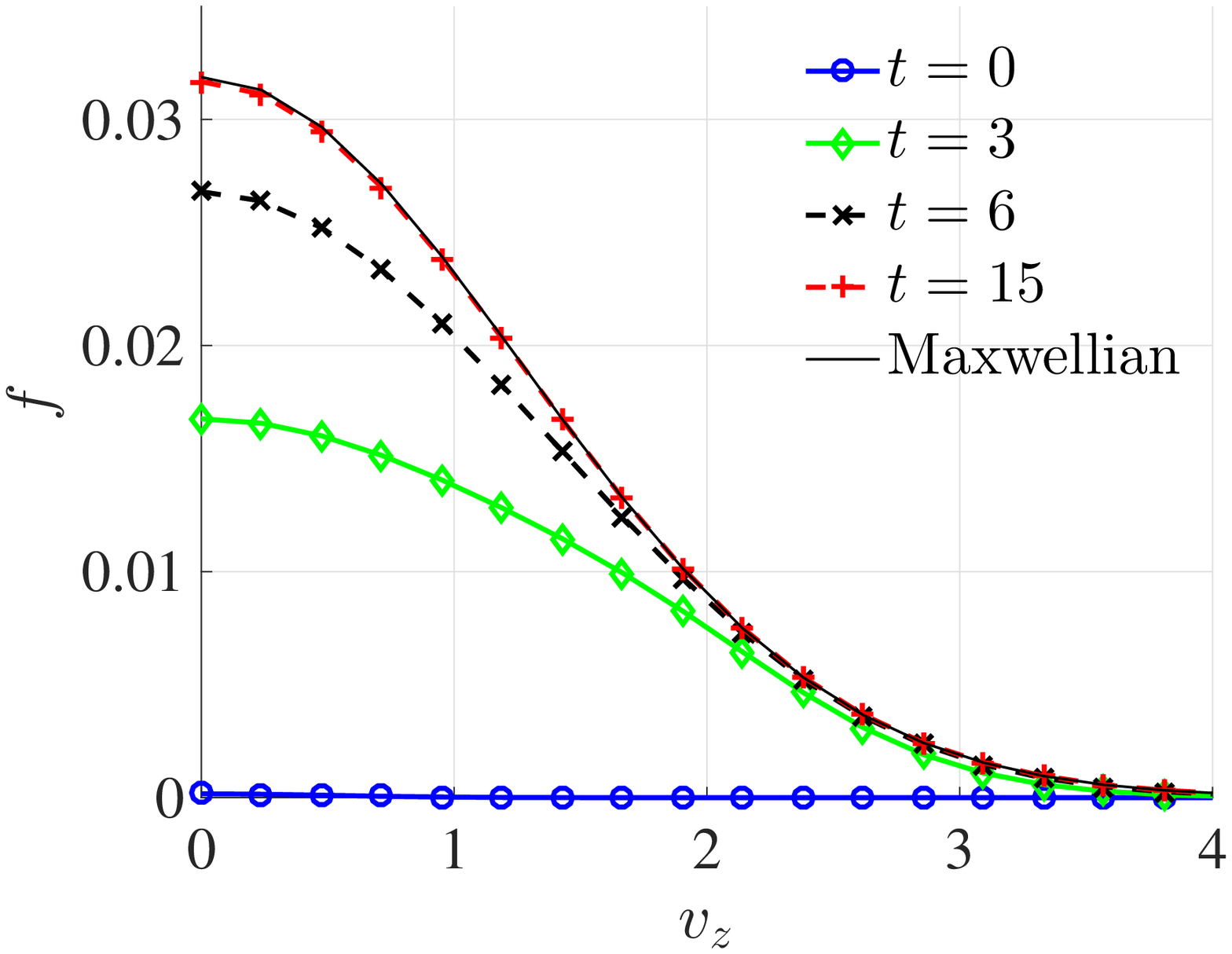}
\includegraphics[width=0.49\textwidth]{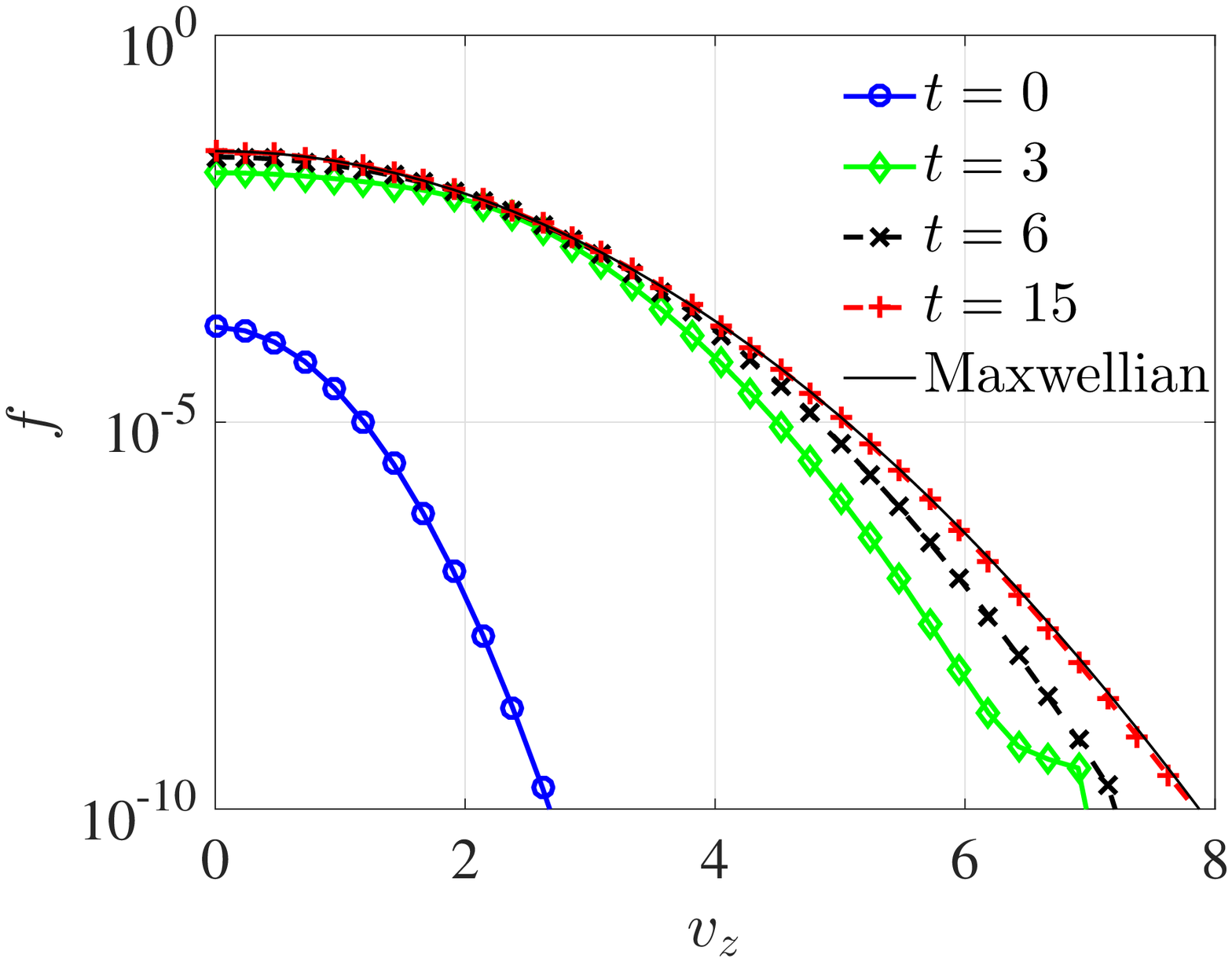}\\
\caption{Velocity pdf on a linear scale (left column) and logarithmic scale (right column),
plotted as a function of $v_x$ (top row) and $v_z$ (bottom row) at the
times shown in the legends for the initial condition in Fig.~\ref{fig:MixMaxX2Contour}.
The limiting Maxwellian pdf is shown with the thin solid black curve. 
%Mix X2 (Job Z4). $g_{\rm tr}=10$, $\Delta t=0.2$, $N=84$.
%These pdfs are all even.  
}
 \label{fig:MixMaxX2TimeSeries}
\end{figure}

In Fig.~\ref{fig:MixMaxX2TimeSeries}, we show the evolution of the
velocity pdf on a linear scale (left column) and logarithmic scale (right column),
plotted as a function of $v_x$ (top row) and $v_z$ (bottom row). 
These results were obtained with $N=84$ and $\Delta t=0.2$. 
We verified our results by comparison to analytic formulae for the 
moments of the velocity pdf~\cite{Rjasanow:2005}. The relative errors in the pressure and 
scalar fourth-order moment were less than $2\times 10^{-4}$ and the absolute error in the
heat flux was less than $4\times 10^{-6}$.
 From $t=0$ to $t=3$, we observe a rapid increase in the very low initial probability of high speed
particles in the $v_z$-direction. Over the same
time period, there is a substantial decay in the  
peaks of  the initial pdf at $\boldsymbol{v}=\boldsymbol{v}_1$ and 
$\boldsymbol{v}=\boldsymbol{v}_2$.
At $t=15$ the agreement with the limiting Maxwellian pdf is excellent 
down to a probability level of $10^{-10}$, i.e., $v < 8$.
However, on a logarithmic scale, when $v > 6$
we observe what appear to be  numerical artifacts 
in the velocity pdf at $t=3$.

\begin{figure}[t!]
\centering
\includegraphics[width=0.49\textwidth]{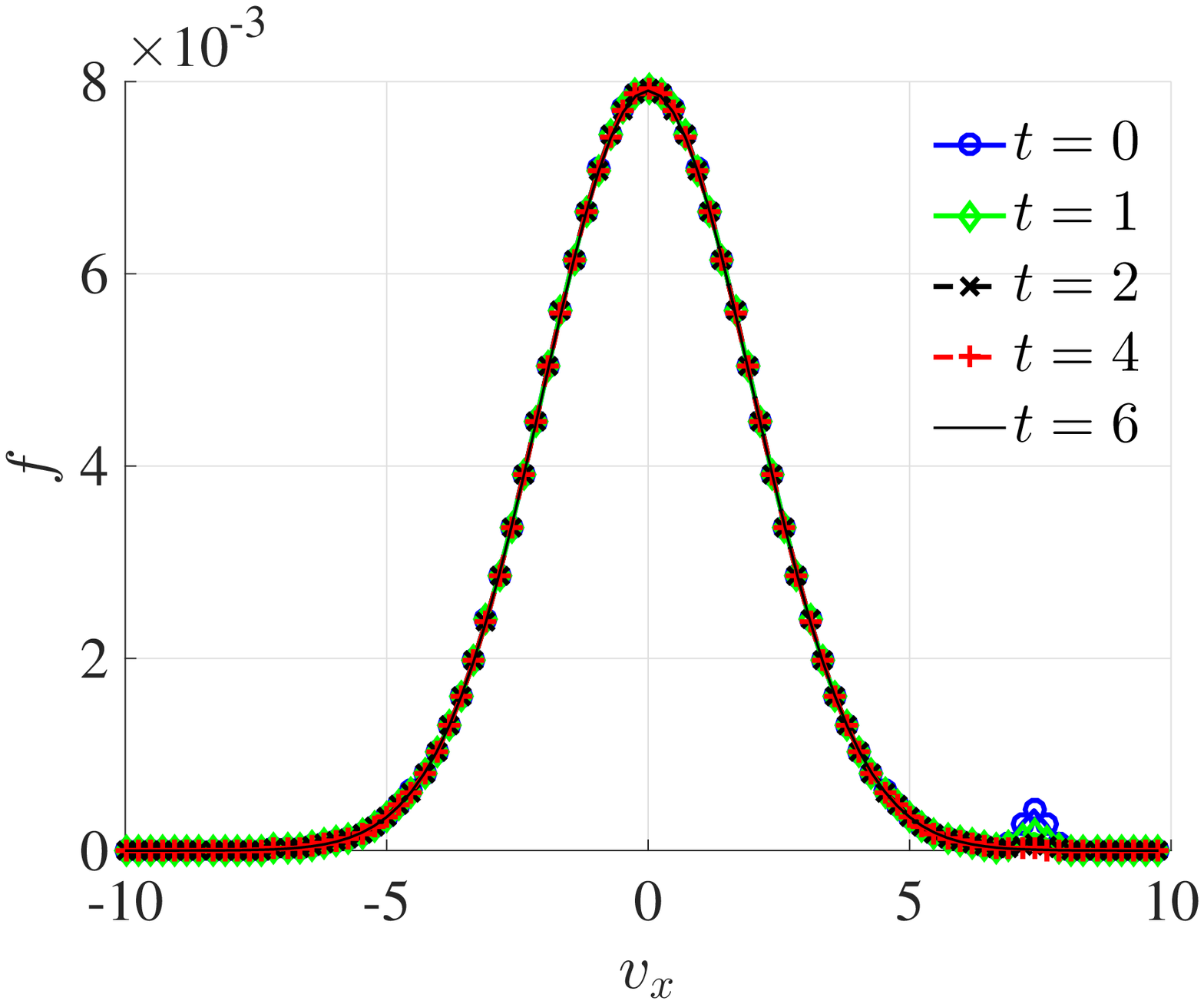}
\includegraphics[width=0.49\textwidth]{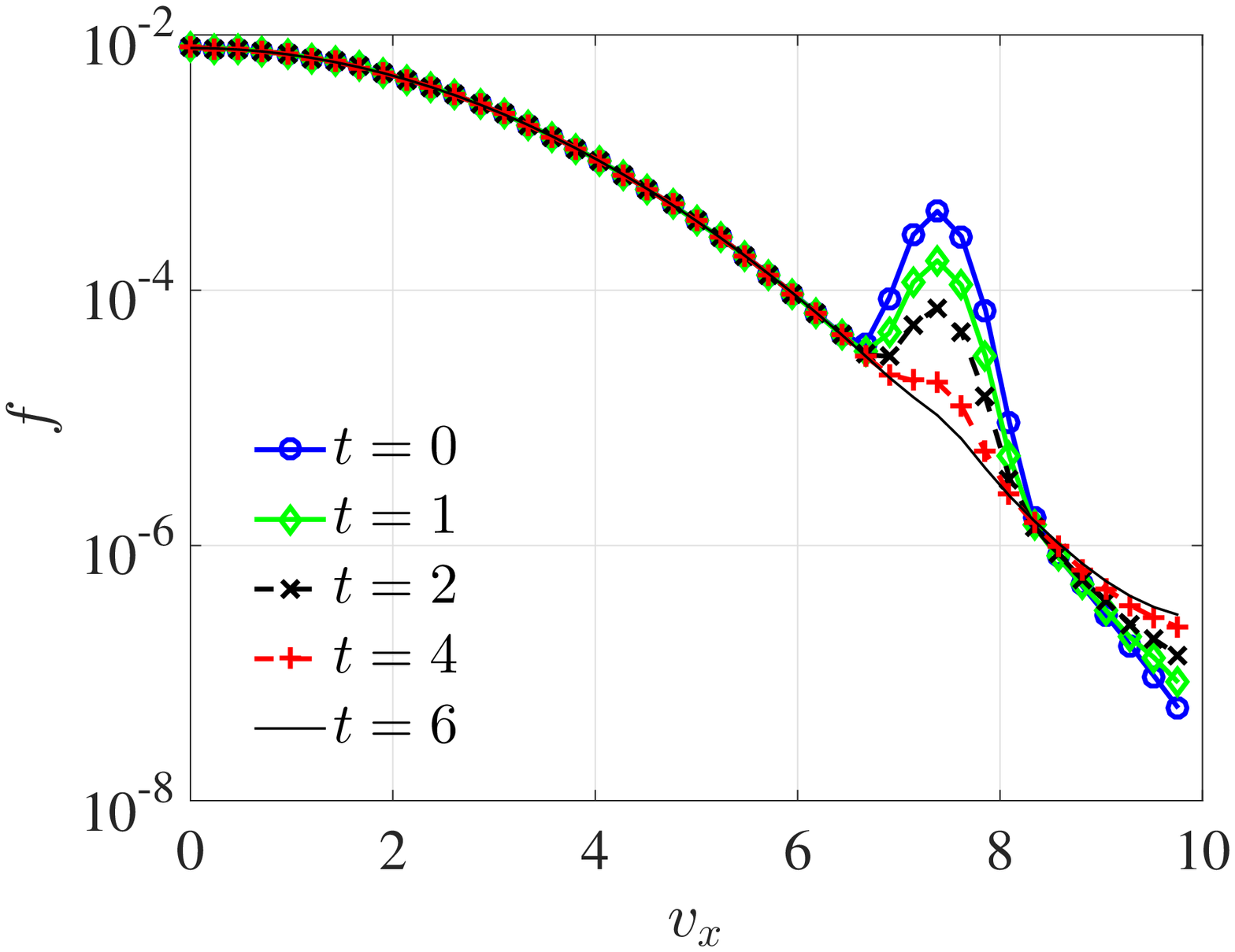}\\
\caption{Velocity pdf on a linear scale (left) and logarithmic scale (right),
plotted as a function of $v_x$  at the
times shown in the legends. 
The parameters in~\eqref{eq:MixMax} were chosen to be
$\omega=0.9999$, $\boldsymbol{v}_1 = (0,0,0)$,
$\boldsymbol{v}_2=(7.38,0,0)$, $T_1=4$, and $T_2=0.0625$. 
%%Mix X9 (Job AD4). $g_{\rm tr}=10$, $\Delta t=0.2$, $N=84$.
}
 \label{fig:MixMaxX10TimeSeries}
\end{figure}

For our second example, we chose $\omega=0.9999$, $\boldsymbol{v}_1 = (0,0,0)$,
$\boldsymbol{v}_2=(7.38,0,0)$, $T_1=4$ and $T_2=0.0625$. With these parameters, the initial pdf is a  perturbation 
of a Maxwellian pdf  which has a small bump centered at 
$\boldsymbol{v}=\boldsymbol{v}_2$ whose amplitude is 0.05 
of that of  the dominant Maxwellian, and which
is located where the dominant Maxwellian has 
a probability density of $10^{-5}$.
Since the probability mass of the bump is negligible, we 
used the dominant Maxwellian rather than the upper bound
of Method~I to estimate the relative error in the truncation of the collision operator.
The resulting contour plot (not pictured) shows that we should choose
$g_{\text{tr}}=10$ to ensure that $\mathcal E_{\operatorname{rel}} < 10^{-1}$ for $v\leq 8$.
As in the previous simulation, we also chose $N=84$ and $\Delta t=0.2$.
The relative error in the pressure was less than $4\times 10^{-4}$
while that of the scalar fourth-order moment was less than $6\times 10^{-2}$.
The absolute error in the heat flux was less than $5\times 10^{-2}$.
In Fig.~\ref{fig:MixMaxX10TimeSeries}, we plot the
velocity pdf at the times shown in the legends 
as a function of $v_x$  for $(v_y,v_z)=(0,0)$
on a linear scale (left) and logarithmic scale (right). 
This simulation result shows the rate at which  
this localized, low-amplitude  perturbation of a Maxwellian pdf relaxes back to the
limiting Maxwellian.  
Because of  how we chose $g_{\text{tr}}$, 
the gradual growth of the pdf where $v> 8.5$ is likely due to 
errors in the numerical computation of the collision operator. 
At all times, the slices of the pdf at $(v_x,v_z)=(0,0)$ and $(v_x,v_y)=(0,0)$ (not shown)
are visually indistinguishable from the dominant Maxwellian.

\subsection{Results for a simple plasma model}

For our final example, we consider a spatially homogeneous
model for the velocity pdf of the electrons in a simplified
plasma system that includes an 
electron gun source, electron-electron collisions, and loss of high velocity electrons into a  wall. We model this system using the equation
\begin{equation}
\frac{\partial f}{\partial t}(t,\mathbf v) \,\,=\,\, Q(f,f)(t,\mathbf v) \,\,+\,\, c_S S(\mathbf v)
 \,\,-\,\, c_L L(\mathbf v)f(t,\mathbf v),
 \label{eq:plasma}
\end{equation}
where the electron gun source is modeled by
$S(\mathbf v) = \exp(\|\mathbf v - \mathbf v_S\|^2/2\sigma_S^2)$
with $\mathbf v_S=(2,0,0)$ and $\sigma_S=0.25$, and the loss 
is given by $L(\mathbf v) = \tfrac{-1}\pi \arctan[ (v_x - v_L)/\sigma_L] + \tfrac 12$
with $v_L = -2$ and $\sigma_L = 10^{-6}$. 
This loss function models absorption of 
particles moving at high speed towards a wall parallel to the $yz$-plane.
To approximately balance gain and loss, we chose the coefficients in \eqref{eq:plasma} to be $c_S=0.1$ and $c_L=10$. For these simulations we chose $g_{\text{tr}}=10$,
$N=80$, and $\Delta t = 0.02$. 
The small value of $\Delta t $ was chosen  to ensure that the numerical solution
did not become negative due to the presence of the loss term.  
The initial velocity pdf was chosen to be a Maxwellian with temperature $T=1$.

\begin{figure}[t!]
\centering
\includegraphics[width=0.49\textwidth]{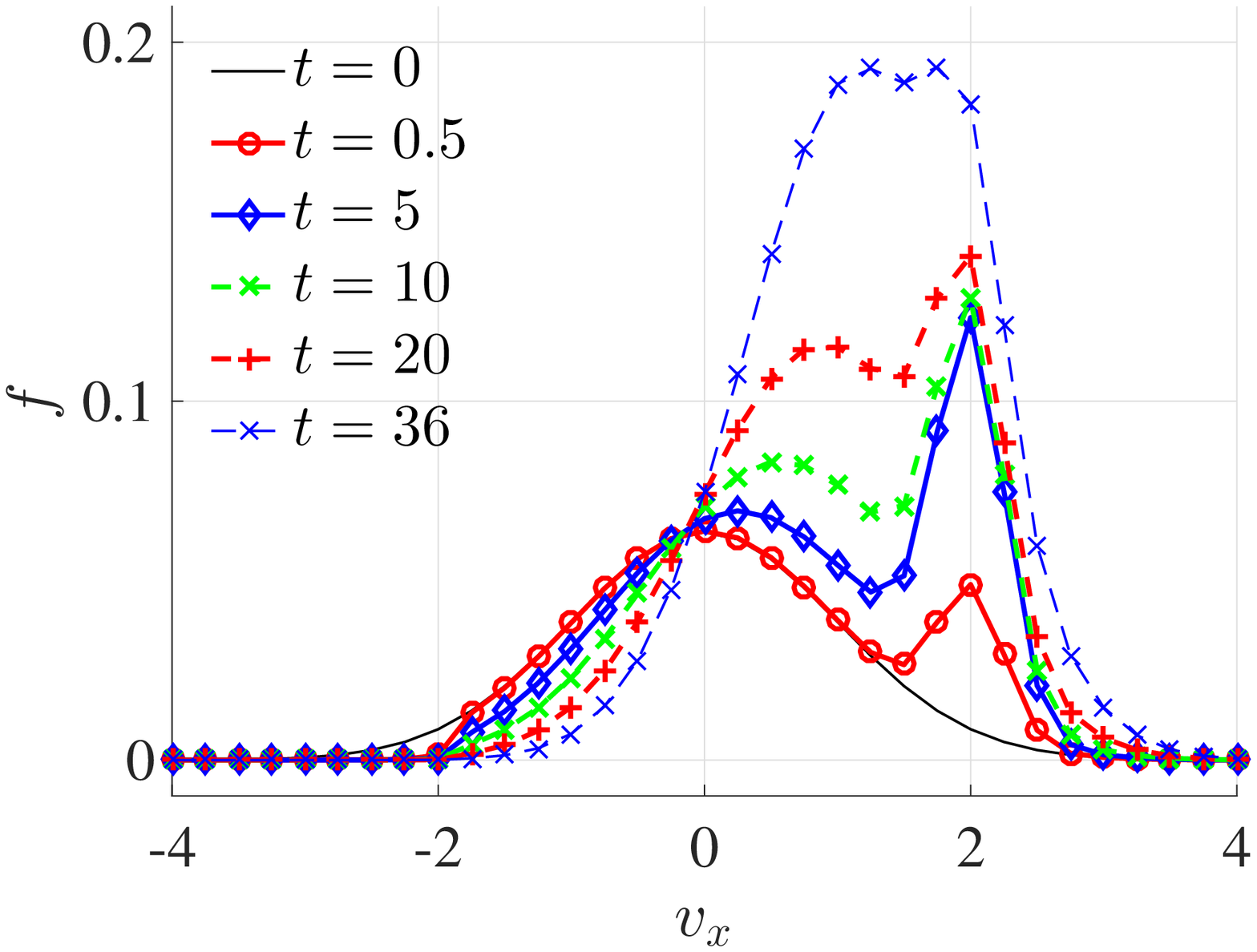}
\includegraphics[width=0.49\textwidth]{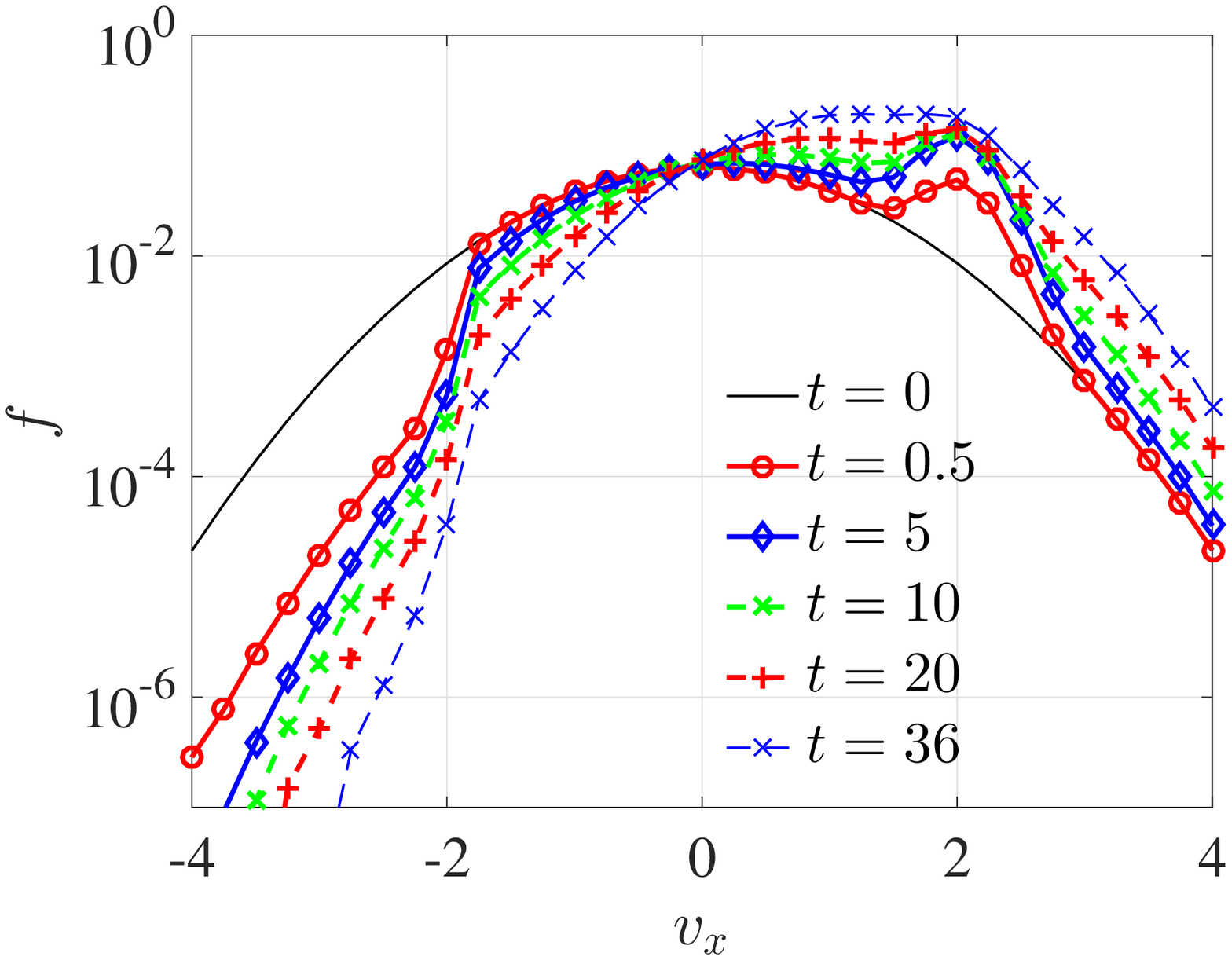}\\
\caption{Velocity pdf for the simple plasma system modeled by \eqref{eq:plasma}
at the times shown in the legends. The pdf is 
plotted as a function of $v_x$  at $(v_y,v_z)=(0,0)$
on a linear scale (left) and a logarithmic scale (right).
%% PlotSourceLossBGSeries.m 
}
 \label{fig:PlasmaTimeSeries}
\end{figure}

In Fig.~\ref{fig:PlasmaTimeSeries}, we plot the time evolution of the velocity pdf
as a function of $v_x$. As time increases from $t=0$ to $t=36$, 
the number density, energy, and the $x$-component of the
momentum all increase due to the source, and the tail of the pdf
in the negative $v_x$-direction deviates significantly from that of a Maxwellian
distribution due to the loss term.  
In addition, the pdf is highly asymmetric in the $v_x$-dimension
due to the combined effects of the source and loss terms.

\section{Conclusions}\label{Sec:Conc}

We have demonstrated the feasibility of using  the spectral-Lagrangian method of
Gamba and Tharkabhushanam  to 
compute the velocity pdf of a particle species  well into the
 low-probability tails.  Calculation of the high-energy tails out to at least three standard deviations could enable improvements to be made in the modeling of chemical reactions 
and ionization events in low-temperature, industrial plasmas.
Although other researchers~\cite{fonn2014hyperbolic,gamba2018galerkin} 
have reported low $L^2$-errors in 
numerical computation of the effect that particle collisions
have on the  distribution of particle velocities, the results presented here are
the first we know of that explicitly study the accuracy of the deterministic computation
of the low-probability tails. 

To obtain these results, we  examined the critical role that the truncation
parameter, $g_{\text{tr}}$, plays in the accuracy of the numerical computation
of the collision operator. Although there is a theoretical guarantee
that the truncated collision operator, $Q^{\text{tr}}$, converges to $Q$ 
as $g_{\text{tr}}\to\infty$,
this result is based on the assumption that the weighted convolution
integral defining $Q^{\text{tr}}$ can be computed exactly without numerical error. 
However, we demonstrated that if $g_{\text{tr}}$ is too large then accurate numerical computation of the  weighted convolution integral  is not feasible 
since the decay rate and degree of oscillation of the convolution weighting function 
both increase as $g_{\text{tr}}$ increases. 
As a consequence, in practice we are forced to examine the trade off between
the error inherent in the truncation of the collision operator and the 
error in the numerical computation of the truncated operator. To do so, 
we derived an upper bound on the pointwise
error between $Q^{\text{tr}}$ and $Q$,
assuming that both operators are computed exactly. 
Unlike in the previous formula for 
$g_{\text{tr}}$ given by 
Gamba and Tharkabhushanam~\cite{GambaThark228p2012}, 
 to obtain this bound we only needed to assume that
 the velocity pdf  is bounded above by a Maxwellian pdf,
 rather than being compactly supported. 
We then showed how to use this bound  to guide the choice of $g_{\text{tr}}$
in numerical computations of the low-probability tails of the velocity pdf.  
Finally,  although our numerical results
were obtained in the spatially homogeneous case, the  error estimate we derived
could also be used to guide the choice of $g_{\text{tr}}$ for the computation
of  spatially inhomogeneous velocity pdfs, 
since the  collision operator is independently  computed at each spatial
position, and, if necessary, the truncation parameter, $g_{\text{tr}}$, can be
chosen to be spatially dependent.

\vskip5pt
\noindent
{\bf Acknowledgements.} We thank Jeff Haack for helpful conversations.

\bibliographystyle{elsarticle-num}
\bibliography{Plasma}
\end{document}